\documentclass[aos,MSNbibl,seceqn,citesort,dvips]{arximspdf}
\usepackage{mathbh}

\doi{10.1214/11-AOS931}
\volume{39}
\issue{6}
\pubyear{2011}
\firstpage{3152}
\lastpage{3181}

\makeatletter

\newcommand{\dsum}{\sum}
\newcommand{\tsum}{\sum}
\newcommand{\dint}{\int}
\newcommand{\bbvarphi}{\bolds{\varphi}}
\newcommand{\bbtau}{\bolds{\tau}}
\newcommand{\bbPsi}{\bolds{\Psi}}
\newcommand{\bbLamda}{\bolds{\Lambda}}
\newcommand{\bbTheta}{\bolds{\Theta}}
\newcommand{\bbvarepsilon}{\bolds{\varepsilon}}
\newcommand{\bbdelta}{\bolds{\delta}}
\newcommand{\bbSigma}{\bolds{\Sigma}}
\newcommand{\bbPhi}{\bolds{\Phi}}
\newcommand{\bbvartheta}{\bolds{\vartheta}}

\newtheorem{theorem}{Theorem}[section]
\newtheorem{lemma}{Lemma}

\makeatother

\begin{document}
\begin{frontmatter}

\title{Gaussian pseudo-maximum likelihood estimation of~fractional time
series models}
\runtitle{Gaussian estimation of fractional time
series models}

\begin{aug}
\author[A]{\fnms{Javier} \snm{Hualde}\corref{}\thanksref{t1}\ead[label=e1]{javier.hualde@unavarra.es}}
\and
\author[B]{\fnms{Peter M.} \snm{Robinson}\thanksref{t2}\ead[label=e2]{P.M.Robinson@lse.ac.uk}}
\runauthor{J. Hualde and P. M. Robinson}
\affiliation{Universidad P\'{u}blica de Navarra and London School of
Economics}
\address[A]{Departamento de Econom\'{i}a\\
Universidad P\'{u}blica de Navarra\\
Campus Arrosad\'{i}a\\
31006 Pamplona\\
Spain\\
\printead{e1}}
\address[B]{Department of Economics\\
London School of Economics\\
Houghton Street\\
London WC2A 2AE\\
United Kingdom\\
\printead{e2}} %adresu isvedimo komanda gale!
\end{aug}

\thankstext{t1}{Supported by the Spanish Ministerio de
Ciencia e Innovaci\'{o}n through a Ram\'{o}n y Cajal contract and ref.
ECO2008-02641.}

\thankstext{t2}{Supported by ESRC Grant
RES-062-23-0036 and Spanish Plan Nacional deI$+$D$+$I Grant
SEJ2007-62908/ECON.}

% HISTORY:
\received{\smonth{5} \syear{2010}}
\revised{\smonth{8} \syear{2011}}

% ABSTRACT
%
\begin{abstract}
We consider the estimation of parametric fractional time series models
in which not only is the memory parameter unknown, but one may not know
whether it lies in the stationary/invertible region or the
nonstationary or noninvertible regions. In these circumstances, a proof
of consistency (which is a prerequisite for proving asymptotic
normality) can be difficult owing to nonuniform convergence of the
objective function over a large admissible parameter space. In
particular, this is the case for the conditional sum of squares
estimate, which can be expected to be asymptotically efficient under
Gaussianity. Without the latter assumption, we establish consistency
and asymptotic normality for this estimate in case of a quite general
univariate model. For a multivariate model, we establish asymptotic
normality of a one-step estimate based on an initial
$\sqrt{n}$-consistent estimate.
\end{abstract}

% KEYWORDS
%

\end{frontmatter}

%s1 #&#
\section{Introduction}

Autoregressive moving average (ARMA) models have featured prominently
in the
analysis of time series. The versions initially stressed in the theoretical
literature (e.g., \cite{hannan,walker}) are stationary and
invertible. Following~\cite{box}, unit root nonstationarity has frequently
been incorporated, while
``overdifferenced'' noninvertible processes have also
featured. Stationary ARMA processes automatically have short memory with
``memory parameter,'' denoted $\delta_{0}$,
taking the value zero, implying a huge behavioral gap relative to unit
root versions, where $\delta_{0}=1$. This has been bridged by
``fractionally-differenced,'' or long
memory, models, a~leading class being the fractional autoregressive
integrated ARMA (FARIMA). A~FARIMA $(p_{1},\delta_{0},p_{2})$ process $
x_{t}$ is given by
%
%e1.1 #&#
%e1.2 #&#
\begin{eqnarray}
\label{a}
x_{t} &=&\Delta^{-\delta_{0}}\{ u_{t}\mathbh{1}(t>0)\} ,\qquad t=0,\pm1,\ldots,
\\
\label{aaaa}
\alpha(L)u_{t} &=&\beta(L)\varepsilon_{t},\qquad t=0,\pm1,\ldots,
\end{eqnarray}
where $\{ x_{t}\} $ is the observable series; $L$ is the lag
operator; $\Delta=1-L$;%
\[
(1-L)^{-\zeta}=\dsum_{j=0}^{\infty}a_{j}(\zeta)L^{j},\qquad
a_{j}(\zeta)=\frac{\Gamma(j+\zeta)}{\Gamma(\zeta)\Gamma(j+1)}
\]
with $\Gamma(\zeta)=\infty$ for $\zeta=0,-1,\ldots,$ and by
convention $%
\Gamma(0)/\Gamma(0)=1$; $\mathbh{1}(\cdot)$ is the indicator function;
$\alpha(L)$ and $\beta (L)$ are real polynomials of degrees $p_{1}$ and
$p_{2}$, which share no common zeros, and all of their zeros are
outside the unit circle in the complex plane; and the $\varepsilon_{t}$
are serially uncorrelated and homoscedastic with zero mean. The reason
(\ref{a}) features the truncated process $u_{t}\mathbh{1} (t>0)$ rather
than simply $u_{t}$ is to simultaneously cover $\delta_{0}$ falling in
both the stationary region $( \delta_{0}< {\frac12} ) $ and the
nonstationary region (\mbox{$\delta_{0}\geq{\frac12}$}, where otherwise the
process would ``blow up''). In the former case, the truncation implies
that $x_{t} $ is only ``asymptotically stationary.'' In recent years,
fractional modeling has found many applications in the sciences and
social sciences; for example, with respect to environmental and
financial data.

Early work on asymptotic statistical theory for fractional models
assumed $ \delta_{0}< {\frac12} $ [and replaced $u_{t}\mathbh{1}(t>0)$
by $u_{t}$ in (\ref{a})]. Assuming $\delta_{0}\in(0, {\frac12} )$,
\cite{dahlhaus,fox,giraitis} and \cite{hosoya1} showed consistency and
asymptotic normality of Whittle estimates (of $\delta_{0}$ and other
parameters, such as the coefficients of $\alpha$ and $\beta$), thereby
achieving analogous results to those of \cite{hannan,walker} for
stationary ARMA processes [i.e., (\ref{aaaa}) with $u_{t}=x_{t}$] and
other short memory models. More recently, \cite{nordman} considered
empirical maximum likelihood inference covering this setting. Note that
\cite{dahlhaus,fox,giraitis} and \cite{hosoya1}, and much other work,
not only excluded $\delta_{0}\geq {\frac12} $ but also the short-memory
case $\delta_{0}=0$, as well as negatively dependent processes where
$\delta_{0}<0$. To some degree, other $\delta_{0}$ can be covered, for
example, for $\delta_{0}\in( 1,{\small3/2}) $ one can first-difference
the data, apply the methods and theory of \cite{dahlhaus,fox,giraitis}
and \cite{hosoya1}, and then add 1 to the memory parameter estimate,
but this still requires prior knowledge that $\delta_{0}$ lies in an
interval of length no more than
$\frac12$.

On the other hand, \cite{beran} argued that the same desirable
properties should hold without so restricting $\delta_{0}$, in case of
a conditional-sum-of-squares estimate, and this would be consistent
with the classical asymptotic properties established by
\cite{robinson1a} for score tests for a unit root and other hypotheses
against fractional alternatives, by comparison with the nonstandard
behavior of unit root tests against stationary autoregressive
alternatives. However, the proof of asymptotic normality in
\cite{beran} appears to assume that the estimate lies in a~small
neighborhood of $\delta_{0}$, without first proving consistency (see
also \cite{tanaka}). Due to a lack\vadjust{\goodbreak} of uniform convergence, consistency
of this implicitly-defined estimate is especially difficult to
establish when the set of admissible values of $\delta$ is large. In
particular, this is the case when $\delta_{0}$ is known only to lie in
an interval of length greater than $\frac12$. In the present paper, we
establish consistency and asymptotic normality when the interval is
arbitrarily large, including (simultaneously) stationary,
nonstationary, invertible and noninvertible values of $\delta _{0}$.
Thus, prior knowledge of which of these phenomena obtains is
unnecessary, and this seems especially practically desirable given, for
example, that estimates near the $\delta_{0}= {\frac12} $ or
$\delta_{0}=1$ boundaries frequently occur in practice, while empirical
interest in autoregressive models with two unit roots suggests
allowance for values in the region of $\delta_{0}=2$ also, and
(following \cite{adensted}) antipersistence and the possibility of
overdifferencing imply the possibility that $\delta_{0}<0$.

We in fact consider a more general model than (\ref{a}), (\ref{aaaa}),
retaining (\ref{a}) but generalizing (\ref{aaaa}) to%
%
%e1.3 #&#
\begin{equation}\label{b}
u_{t}=\theta(L;\bbvarphi_{0})\varepsilon_{t},\qquad t=0,\pm1,\ldots,
\end{equation}
where $\varepsilon_{t}$ is a zero-mean unobservable white noise
sequence, $%
\bbvarphi_{0}$ is an unknown $p\times1$ vector, $\theta
(s;\bbvarphi)=\dsum_{j=0}^{\infty}\theta_{j}(\bbvarphi)s^{j}$%
, where for all $\bbvarphi$, $\theta_{0}(\bbvarphi)=1$, $%
\theta(s;\bbvarphi)\dvtx\mathbb{\mathbb{C}
}\times\mathbb{R}^{p}$ is continuous in $s$ and $\vert\theta(s;%
\bbvarphi)\vert\neq0$, $\vert s\vert\leq1$.
More detailed conditions will be imposed below. The role of $\theta$
in (%
\ref{b}), like $\alpha$ and $\beta$ in~(\ref{aaaa}), is to permit
parametric short memory autocorrelation. We allow for the simplest case
FARIMA$(0,\delta_{0},0)$ by taking $\bbvarphi_{0}$ to be empty.
Another model covered by~(\ref{b}) is the exponential-spectrum one of
\cite%
{bloomfield} (which in conjunction with fractional differencing leads
to a
relatively neat covariance matrix formula~\cite{robinson1a}). Semiparametric
models (where $u_{t}$ has nonparametric autocovariance structure; see, e.g.,
\cite{robinson1b,shimotsu}) afford still greater flexibility than~(%
\ref{b}), but also require larger samples in order for comparable precision
to be achieved. In more moderate-sized samples, investment in a parametric
model can prove worthwhile, even the simple FARIMA($1$, $\delta_{0}$, $0$)
employed in the Monte Carlo simulations reported in the supplementary
material~\cite{hualde}, while model choice procedures can be employed to
choose $p_{1}$ and $p_{2}$ in the FARIMA($p_{1}, \delta_{0},p_{2}$), as
illustrated in the empirical examples included in the supplementary material
\cite{hualde}.

We wish to estimate $\bbtau_{0}=(\delta_{0},\bbvarphi%
_{0}^{\prime})^{\prime}$ from observations $x_{t}$, $t=1,\ldots,n$. For any
admissible $\bbtau=(\delta,\bbvarphi^{\prime
})^{\prime}$,
define
%
%e1.4 #&#
\begin{equation}\label{d}
\varepsilon_{t}(\bbtau)=\Delta^{\delta}\theta^{-1}(L;\bbvarphi)x_{t},\qquad t\geq1,
\end{equation}
noting that (\ref{a}) implies $x_{t}=0$, $t\leq0$. For a given user-chosen
optimizing set $\mathcal{T}$, define as an estimate of $\bbtau_{0}$%
%
%e1.5 #&#
\begin{equation} \label{e}
\widehat{\bbtau}=\mathop{\arg\min}_{\bbtau\in\mathcal{T}}
R_{n}(\bbtau),
\end{equation}
where
%
%e1.6 #&#
\begin{equation} \label{f}
R_{n}(\bbtau)=\frac{1}{n}\tsum_{t=1}^{n}\varepsilon_{t}^{2}(%
\bbtau),\vadjust{\goodbreak}
\end{equation}
and $\mathcal{T}=\mathcal{I}\times\Psi$, where $\mathcal{I}=\{ \delta
\dvtx\bigtriangledown_{1}\leq\delta\leq\bigtriangledown_{2}\} $ for
given $\bigtriangledown_{1}$, $\bigtriangledown_{2}$ such that $%
\bigtriangledown_{1}<\bigtriangledown_{2}$, $\Psi$ is a compact
subset of
$\mathbb{R}^{p}$ and $\bbtau_{0}\in\mathcal{T}$.

The estimate $\widehat{\bbtau}$ is sometimes termed
``conditional sum of squares'' (though
``truncated sum of squares'' might be more
suitable). It has the anticipated advantage of having the same limit
distribution as the maximum likelihood estimate of $\bbtau_{0}$
under Gaussianity, in which case it is asymptotically efficient (though here
we do not assume Gaussianity). It was employed by \cite{box} in estimation
of nonfractional ARMA models (when $\delta_{0}$ is a given integer), by
\cite{li,robinson3} in stationary FARIMA models, where $0<\delta
_{0}<1/2$, and by \cite{beran,tanaka} in nonstationary FARIMA
models, allowing $\delta_{0}\geq1/2$.

The following section sets down detailed regularity conditions, a formal
statement of asymptotic properties and the main proof details. Section \ref{sec3}
provides asymptotically normal estimates in a multivariate extension of~(\ref{a}),~(\ref{b}). Joint modeling of related processes is important both for
reasons of parsimony and interpretation, and multivariate fractional
processes are currently relatively untreated, even in the stationary case.
Further possible extensions are discussed in Section \ref{sec4}. Useful lemmas are
stated in Section \ref{sec5}. Due to space restrictions, the proofs of these lemmas,
along with an analysis of finite-sample performance of the procedure
and an
empirical application, are included in the supplementary material
\cite{hualde}.

%s2 #&#
\section{Consistency and asymptotic normality}\label{sec2}

%s2.1 #&#
\subsection{\texorpdfstring{Consistency of $\widehat{\bbtau}$}{Consistency of tau}}

Our first two assumptions will suffice for consistency.

\begin{enumerate}[A1.]
\item[A1.]

\begin{enumerate}[(iii)]
\item[(i)]
\[
\vert\theta( s;\bbvarphi) \vert\neq
\vert\theta( s;\bbvarphi_{0}) \vert
\]
for all $\bbvarphi\neq\bbvarphi_{0}$, $\bbvarphi\in
\Psi$, on a set $S\subset\{ s\dvtx\vert s\vert=1\} $
of positive Lebesgue measure;

\item[(ii)] for all $\bbvarphi$, $\theta( e^{i\lambda};%
\bbvarphi) $ is differentiable in $\lambda$ with derivative
in $\operatorname{Lip}( \varsigma) $, $\varsigma>1/2;$

\item[(iii)] for all $\lambda$, $\theta( e^{i\lambda};\bbvarphi) $ is continuous in $\bbvarphi;$

\item[(iv)] for all $\bbvarphi\in\Psi$, $\vert\theta(
s;\bbvarphi) \vert\neq0, \vert s\vert\leq
1$.
\end{enumerate}
\end{enumerate}

Condition (i) provides identification while (ii) and (iv) ensure that $u_{t}$
is an invertible short-memory process (with spectrum that is bounded and
bounded away from zero at all frequencies). Further, by (ii) the derivative
of~$\theta(e^{i\lambda};\bbvarphi)$ has Fourier coefficients
$j\theta_{j}( \bbvarphi) =O( j^{-\varsigma}) $
as $j\rightarrow\infty$, for all $\bbvarphi$, from page 46 of
\cite{zygmund}, so that, by compactness of $\Psi$ and continuity of
$\theta
_{j}( \bbvarphi) $ in $\bbvarphi$ for all $j$,
%
%e2.1 #&#
\begin{equation} \label{h14}
{\sup_{\bbvarphi\in\Psi}}\vert\theta_{j}( \bbvarphi) \vert=O\bigl( j^{-(
1+\varsigma) }\bigr) \qquad\mbox{as }j\rightarrow\infty.
\end{equation}
Also, writing $\theta^{-1}( s;\bbvarphi) =\phi( s;%
\bbvarphi) =\dsum_{j=0}^{\infty}\phi_{j}(
\bbvarphi) s^{j}$, we have $\phi_{0}( \bbvarphi%
) =1$ for all~$\bbvarphi$, and (ii), (iii) and (iv) imply that
%
%e2.2 #&#
\begin{equation} \label{eee}
{\sup_{\bbvarphi\in\Psi}}\vert\phi_{j}( \bbvarphi%
) \vert=O\bigl( j^{-( 1+\varsigma) }\bigr)\qquad \mbox{as
}j\rightarrow\infty.
\end{equation}
Finally, (ii) also implies that%
%
%e2.3 #&#
\begin{equation} \label{zf}
{\mathop{\inf_{\vert s\vert=1}}_{\bbvarphi\in\Psi}}
\vert\phi( s;\bbvarphi) \vert>0.
\end{equation}
Assumption A1 is easily satisfied by standard parameterizations of stationary and
invertible ARMA processes (\ref{aaaa}) in which autoregressive and moving
average orders are not both over-specified. More generally, A1 is
similar to
conditions employed in asymptotic theory for the estimate $\widehat
{\bbtau}$ and other forms of Whittle estimate that restrict to stationarity
(see, e.g., \cite{dahlhaus,fox,giraitis,hosoya1,robinson3}) and not
only is it readily verifiable because $\theta$ is a
known parametric function, but in practice $\theta$ satisfying A1 are
invariably employed by practitioners.

\begin{enumerate}[A2.]
\item[A2.] The $\varepsilon_{t}$ in (\ref{b}) are stationary and ergodic
with finite fourth moment, and%
%
%e2.4 #&#
\begin{equation} \label{28}
E( \varepsilon_{t}\vert\mathcal{F}_{t-1}) =0,\qquad
E( \varepsilon_{t}^{2}\vert\mathcal{F}_{t-1})
=\sigma_{0}^{2}
\end{equation}
almost surely, where $\mathcal{F}_{t}$ is the $\sigma$-field of events
generated by $\varepsilon_{s}$, $s\leq t$, and conditional (on
$\mathcal{F}%
_{t-1}$) third and fourth moments of $\varepsilon_{t}$ equal the
corresponding unconditional moments.
\end{enumerate}

Assumption A2 avoids requiring independence or identity of distribution of
$\varepsilon
_{t}$, but rules out conditional heteroskedasticity. It has become fairly
standard in the time series asymptotics literature since \cite{hannan}.
%
%th2.1 #&#
\begin{theorem}
Let (\ref{a}), (\ref{b}) and \textup{A1, A2} hold. Then as $n\rightarrow\infty$
%
%e2.5 #&#
\begin{equation} \label{zd}
\widehat{\bbtau}\rightarrow_{p}\bbtau_{0}.
\end{equation}
\end{theorem}
\begin{pf}
We give the proof for the most general case where $\bigtriangledown
_{1}<\delta_{0}-\frac{1}{2}$, but our proof trivially covers the $%
\bigtriangledown_{1}\geq\delta_{0}-\frac{1}{2}$ situation, for which some
of the steps described below are superfluous. The proof begins standardly.
For $\varepsilon>0$, define $N_{\varepsilon}=\{ \bbtau%
\dvtx\Vert\bbtau-\bbtau_{0}\Vert<\varepsilon\} $, $%
\overline{N}_{\varepsilon}=\{ \bbtau\dvtx\bbtau\notin
N_{\varepsilon},\bbtau\in\mathcal{T}\} $. For small enough $%
\varepsilon$,%
%
%e2.6 #&#
\begin{equation} \label{new1}
\Pr( \widehat{\bbtau}\in\overline{N}_{\varepsilon}) \leq
\Pr\Bigl( \inf_{\bbtau\in\overline{N}_{\varepsilon}}S_{n}(
\bbtau) \leq0\Bigr) ,
\end{equation}
where $S_{n}( \bbtau) =R_{n}( \bbtau)
-R_{n}( \bbtau_{0}) $. The remainder of the proof
reflects the fact that $R_{n}(\bbtau)$, and thus $S_{n}(\bbtau)$, converges in probability to a well-behaved function when $\delta
>\delta_{0}-{\frac12}$, and diverges when $\delta<\delta_{0}- {\frac12}$, while the
need to establish uniform convergence, especially in a neighborhood
of
$\delta=\delta_{0}- {\frac12}$, requires additional special treatment.
Consequently,\vspace*{1pt} for arbitrarily\vadjust{\goodbreak} small \mbox{$\eta>0$}, such that
$\eta<\delta_{0}-\frac {1}{2}-\bigtriangledown_{1}$, we define the
nonintersecting sets $\mathcal{I}_{1}=\{
\delta\dvtx\bigtriangledown_{1}\leq\delta\leq
\delta_{0}-\frac{1}{2}-\eta\} $, $\mathcal{I}_{2}=\{ \delta\dvtx\delta
_{0}-\frac{1}{2%
}-\eta<\delta<\delta_{0}-\frac{1}{2}\} $, $\mathcal{I}_{3}=\{
\delta\dvtx\delta_{0}-\frac{1}{2}\leq\delta\leq\delta_{0}-\frac
{1}{2}+\eta
\} $, $\mathcal{I}_{4}=\{ \delta\dvtx\delta_{0}-\frac{1}{2}+\eta
<\delta\leq\bigtriangledown_{2}\} $. Correspondingly, define
$\mathcal{T}_{i}=\mathcal{I}_{i}\times
\Psi$, $i=1,\ldots,4$, so $\mathcal{T}=\bigcup_{i=1}^{4}\mathcal{T}_{i} $. Thus,
from (\ref{new1}) it remains to prove
%
%e2.7 #&#
\begin{equation} \label{new2}
\Pr\Bigl( \inf_{\bbtau\in\overline{N}_{\varepsilon}\cap\mathcal
{T}%
_{i}}S_{n}( \bbtau) \leq0\Bigr) \rightarrow0\qquad\mbox{as
}n\rightarrow\infty,\qquad i=1,\ldots,4.
\end{equation}
Each of the four proofs differs, and we describe them in reverse order.

\textit{Proof of} (\ref{new2}) \textit{for $i=4$.} By a familiar
argument, the result follows if for $\bbtau\in\mathcal{T}_{4}$
there is a deterministic function $U( \bbtau) $ (not
depending on $n$), such that
\[
S_{n}( \bbtau) =U( \bbtau)
-T_{n}( \bbtau) ,
\]
where
%
%e2.8 #&#
\begin{equation} \label{1}
\inf_{\overline{N}_{\varepsilon}\cap\mathcal{T}_{4}}U( \bbtau%
) >\epsilon,
\end{equation}
$\epsilon$ throughout denoting a generic arbitrarily small positive
constant, and
%
%e2.9 #&#
\begin{equation} \label{2}
{\sup_{\mathcal{T}_{4}}}\vert T_{n}( \bbtau)
\vert=o_{p}( 1) .
\end{equation}
Since $x_{t}=0$, $t\leq0$, for $\bbtau\in\mathcal{T}_{4}$
we set [cf. (\ref{d})], $\zeta_{t}( \bbtau) =\Delta
^{\delta-\delta_{0}}\phi( L;\bbvarphi) u_{t}$, $%
U( \bbtau) =E\zeta_{t}^{2}( \bbtau)
-\sigma_{0}^{2}$ and $T_{n}( \bbtau) =R_{n}(
\bbtau_{0}) -\sigma_{0}^{2}-\{ R_{n}( \bbtau%
) -E\zeta_{t}^{2}( \bbtau) \} $. We may
write
\[
U( \bbtau) =\sigma_{0}^{2}\biggl( \frac{1}{2\pi}%
\int_{-\pi}^{\pi}\frac{g( \lambda) }{g_{0}(
\lambda) }\,d\lambda-1\biggr) ,
\]
where
\[
g( \lambda) =\vert1-e^{i\lambda}\vert^{2(
\delta-\delta_{0}) }\vert\phi( e^{i\lambda};\bbvarphi) \vert^{2},\qquad
g_{0}( \lambda)= g( \lambda) \vert_{\bbtau=\bbtau%
_{0}}.
\]
For all $\bbtau$ $( 2\pi) ^{-1}\dint_{-\pi}^{\pi}\log( g(
\lambda) /g_{0}( \lambda) ) \,d\lambda=0$, so by Jensen's inequality
%
%e2.10 #&#
\begin{equation} \label{ze}
\frac{1}{2\pi}\int_{-\pi}^{\pi}\frac{g( \lambda) }{%
g_{0}( \lambda) }\,d\lambda\geq1.
\end{equation}
Under A1(i), we have strict inequality in (\ref{ze}) for all $\bbtau%
\neq\bbtau_{0}$, so that by continuity in $\bbtau$ of the
left-hand side of (\ref{ze}), (\ref{1}) holds. Next, write
\[
\varepsilon_{t}( \bbtau)
=\sum_{j=0}^{t-1}c_{j}( \bbtau) u_{t-j},\qquad
\zeta_{t}( \bbtau) =\sum_{j=0}^{\infty
}c_{j}( \bbtau) u_{t-j},
\]
where $c_{j}( \bbtau) =\dsum_{k=0}^{j}\phi
_{k}( \bbvarphi) a_{j-k}( \delta_{0}-\delta)
$. Because, given A2, the $\varepsilon_{t}^{2}-\sigma_{0}^{2}$ are
stationary martingale
differences,
%
%e2.11 #&#
\begin{equation} \label{1010}
R_{n}( \bbtau_{0}) -\sigma_{0}^{2}=\frac{1}{n}%
\sum_{t=1}^{n}( \varepsilon_{t}^{2}-\sigma_{0}^{2})
\rightarrow_{p}0\qquad\mbox{as }n\rightarrow\infty.
\end{equation}
Then defining $\gamma_{k}=E( u_{t}u_{t-k}) $, and henceforth
writing $c_{j}=c_{j}( \bbtau) $, (\ref{2}) would hold on
showing that
%
%e2.12 #&#
%e2.13 #&#
%e2.14 #&#
\begin{eqnarray}
\label{4}
\sup_{\mathcal{T}_{4}}\Biggl\vert\frac{1}{n}\sum_{t=1}^{n}\Biggl[
\Biggl( \sum_{j=0}^{t-1}c_{j}u_{t-j}\Biggr) ^{2}-E\Biggl(
\sum_{j=0}^{t-1}c_{j}u_{t-j}\Biggr) ^{2}\Biggr] \Biggr\vert
&=&o_{p}( 1) ,\\[-2pt]
\label{5}
\sup_{\mathcal{T}_{4}}\Biggl\vert\frac{1}{n}\sum_{t=1}^{n}\sum%
_{j=0}^{t-1}\sum_{k=t}^{\infty}c_{j}c_{k}\gamma
_{j-k}\Biggr\vert&=&o_{p}( 1) , \\[-2pt]
\label{5bis}
\sup_{\mathcal{T}_{4}}\Biggl\vert\frac{1}{n}\sum_{t=1}^{n}\sum%
_{j=t}^{\infty}\sum_{k=t}^{\infty}c_{j}c_{k}\gamma
_{j-k}\Biggr\vert&=&o_{p}( 1) .
\end{eqnarray}
We first deal with (\ref{4}). The term whose modulus is taken is
%
%e2.15 #&#
\begin{eqnarray} \label{ai}
&&\frac{1}{n}\sum_{j=0}^{n-1}c_{j}^{2}\sum_{l=1}^{n-j}(
u_{l}^{2}-\gamma_{0})\nonumber\\[-2pt]
&&\quad{} +\frac{2}{n}\sum_{j=0}^{n-2}\sum%
_{k=j+1}^{n-1}c_{j}c_{k}\sum_{l=k-j+1}^{n-j}\bigl\{
u_{l}u_{l-( k-j) }-\gamma_{j-k}\bigr\}\\[-2pt]
&&\qquad
=( a) +( b) .\nonumber
\end{eqnarray}
First,%
\[
{E\sup_{\mathcal{T}_{4}}}\vert( a) \vert\leq\frac{1}{n%
}\dsum_{j=0}^{n-1}\sup_{\mathcal{T}_{4}}c_{j}^{2}E\Biggl\vert
\dsum_{l=1}^{n-j}( u_{l}^{2}-\gamma_{0}) \Biggr\vert.
\]
It can be readily shown that, uniformly in $j$, $\operatorname{Var}(
\dsum_{l=1}^{n-j}u_{l}^{2}) =O( n) $, so
\[
{\sup_{\mathcal{T}_{4}}}\vert( a) \vert=O_{p}\Biggl(
n^{-{1/2}}\sum_{j=1}^{\infty}j^{-2\eta-1}\Biggr) =O_{p}(
n^{-{1/2}})
\]
by Lemma \ref{lemma1}. Next, by summation by parts, $( b) $ is equal to
\begin{eqnarray*}
&&\frac{2c_{n-1}}{n}\sum_{j=0}^{n-2}c_{j}\sum_{k=j+1}^{n-1}%
\sum_{l=k-j+1}^{n-j}\bigl\{ u_{l}u_{l-( k-j) }-\gamma
_{j-k}\bigr\} \\[-2pt]
&&\quad{}-\frac{2}{n}\sum_{j=0}^{n-2}c_{j}\sum_{k=j+1}^{n-2}(
c_{k+1}-c_{k})
\sum_{r=j+1}^{k}\sum_{l=r-j+1}^{n-j}\bigl\{ u_{l}u_{l-(
r-j) }-\gamma_{j-r}\bigr\} \\[-2pt]
&&\qquad=( b_{1}) +( b_{2}) .
\end{eqnarray*}
It can be easily shown that, uniformly in $j$,%
\[
\operatorname{Var}\Biggl(
\sum_{k=j+1}^{n-1}\sum_{l=k-j+1}^{n-j}u_{l}u_{l-(
k-j) }\Biggr) =O( n^{2}) ,\vadjust{\goodbreak}
\]
so we have%
\begin{eqnarray*}
{E\sup_{\mathcal{T}_{4}}}\vert( b_{1}) \vert&\leq&
Kn^{-\eta-{3/2}}\\
&&{}\times\sum_{j=1}^{n}j^{-\eta-{1/2}}\Biggl\{
\operatorname{Var}\Biggl(
\sum_{k=j+1}^{n-1}\sum_{l=k-j+1}^{n-j}u_{l}u_{l-(
k-j) }\Biggr) \Biggr\} ^{1/2}\\
&\leq& Kn^{-2\eta}
\end{eqnarray*}
by Lemma \ref{lemma1}, where $K$ throughout denotes a generic finite but arbitrarily
large positive constant. Similarly,
\begin{eqnarray*}
&&{E\sup_{\mathcal{T}_{4}}}\vert( b_{2}) \vert\\
&&\qquad\leq
Kn^{-1}\sum_{j=1}^{n}j^{-\eta-{1/2}}\sum_{k=j+1}^{n}k^{%
\max( -\eta-{3/2},-( 1+\varsigma) ) }\\
&&\qquad\quad\hspace*{107pt}{}\times\Biggl\{
\operatorname{Var}\Biggl( \sum_{r=j+1}^{k}\sum_{l=r-j+1}^{n-j}u_{l}u_{l-(
r-j) }\Biggr) \Biggr\}^{1/2}
\end{eqnarray*}
by Lemma \ref{lemma1}, where $\varsigma$ was introduced in A1(ii). It can be readily
shown that%
\[
\operatorname{Var}\Biggl( \sum_{r=j+1}^{k}\sum_{l=r-j+1}^{n-j}u_{l}u_{l-(
r-j) }\Biggr) \leq K( k-j) ( n-j) .
\]
Take $\eta$ such that $\eta+\frac{3}{2}<1+\varsigma$. Then%

\begin{eqnarray*}
{E\sup_{\mathcal{T}_{4}}}\vert( b_{2}) \vert &\leq
&Kn^{-{1/2}}\dsum_{j=1}^{n}j^{-\eta-{1/2}%
}\dsum_{k=j+1}^{n}k^{-\eta-{3/2}}( k-j) ^{{1/2%
}} \\
&\leq&Kn^{-{1/2}}\dsum_{j=1}^{n}j^{-\eta-{1/2}%
}\dsum_{k=1}^{n}( k+j) ^{-\eta-{3/2}}k^{{1/2}%
}.
\end{eqnarray*}
This is bounded by%
%
%e2.16 #&#
\begin{equation} \label{new3}
Kn^{-{1/2}}\dsum_{j=1}^{n}j^{-3\eta-{1/2}%
}\dsum_{k=1}^{n}k^{\eta-1},
\end{equation}
because $( k+j) ^{-\eta-{3/2}}\leq j^{-2\eta}k^{\eta-%
{3/2}}$. For small enough $\eta$, (\ref{new3}) is bounded by $%
Kn^{-2\eta}$, to complete the proof of (\ref{4}).
Next, the term whose modulus is taken in (\ref{5}) is
%
%e2.17 #&#
\begin{equation} \label{ab}
\frac{1}{n}\sum_{t=1}^{n}\int_{-\pi}^{\pi}f( \lambda
) \sum_{j=0}^{t-1}\sum_{k=t}^{\infty
}c_{j}c_{k}e^{i( j-k) \lambda}\,d\lambda,
\end{equation}
where $f( \lambda) $ denotes the spectral density of $u_{t}$. By
boundedness of $f$
(implied by assumption A1) and the Cauchy inequality, (\ref{ab}) is
bounded by
\begin{eqnarray*}
&&
Kn^{-1}\sum_{t=1}^{n}\Biggl\{ \int_{-\pi}^{\pi}\Biggl\vert
\sum_{j=0}^{t-1}c_{j}e^{ij\lambda}\Biggr\vert^{2}\,d\lambda
\int_{-\pi}^{\pi}\Biggl\vert\sum_{k=t}^{\infty
}c_{k}e^{-ik\lambda}\Biggr\vert^{2}\,d\lambda\Biggr\}^{{1/2}}\\
&&\qquad\leq
Kn^{-1}\sum_{t=1}^{n}\Biggl\{
\sum_{j=0}^{t-1}c_{j}^{2}\sum_{k=t}^{\infty}c_{k}^{2}\Biggr\}
^{{1/2}},
\end{eqnarray*}
so the left-hand side of (\ref{5}) is bounded by
\[
Kn^{-1}\sum_{t=1}^{n}\Biggl\{
\sum_{j=1}^{t}j^{-2\eta-1}\sum_{k=t}^{\infty
}k^{-2\eta-1}\Biggr\}^{{1/2}}\leq
Kn^{-1}\sum_{t=1}^{n}t^{-\eta}\leq
Kn^{-\eta}=o( 1)
\]
by Lemma \ref{lemma1}, to establish (\ref{5}).
Finally, by a similar reasoning, the term whose modulus is taken in
(\ref%
{5bis}) is bounded by
\[
Kn^{-1}\sum_{t=1}^{n}\Biggl\{ \int_{-\pi}^{\pi}\Biggl\vert
\sum_{j=t}^{\infty}c_{j}e^{ij\lambda}\Biggr\vert^{2}\,d\lambda
\Biggr\}^{{1/2}}\leq Kn^{-1}\sum_{t=1}^{n}t^{-2\eta}\leq Kn^{-2\eta}
\]
to conclude the proof of (\ref{5bis}), and thence of (\ref{2}). Thus,
(\ref{new2}) is proved for
$i=4$.
With respect to (\ref{new2}) for $i=1,2,3$, note from $\mathcal
{T}_{i}\cap
\overline{N}_{\varepsilon}\equiv\mathcal{T}_{i}$ for such $i$, and
(\ref
{1010}), that these results follow if
%
%e2.18 #&#
\begin{equation} \label{ac}
\Pr\Bigl( \inf_{\mathcal{T}_{i}}R_{n}( \bbtau)
\leq K \Bigr) \rightarrow0\qquad\mbox{as }n\rightarrow\infty,
i=1,2,3.
\end{equation}

\textit{Proof of} (\ref{new2}) \textit{for $i=3$}. Denote, for any sequence
$\zeta_{t}$, $w_{\zeta}( \lambda)
=\break n^{-{1/2}}\dsum_{t=1}^{n}\zeta_{t}\times\allowbreak e^{it\lambda}$, $%
I_{\zeta}( \lambda) =\vert w_{\zeta}( \lambda
) \vert^{2}$, the discrete Fourier transform and periodogram,
respectively, and $\lambda
_{j}=2\pi j/n$. For $V_{n}( \bbtau) $ satisfying Lemma \ref{lemma3},
setting $\bbtau^{\ast}=( \delta,\bbvarphi_{0}^{\prime
}) ^{\prime}$,
\[
R_{n}( \bbtau) =\frac{1}{n}\sum_{j=1}^{n}I_{\varepsilon
( \bbtau) }( \lambda_{j})
=\frac{1}{n}\sum_{j=1}^{n}\vert\xi( e^{i\lambda_{j}};\bbvarphi) \vert^{2}I_{\varepsilon(
\bbtau^{\ast
}) }( \lambda_{j}) +\frac{1}{n}V_{n}( \bbtau) ,
\]
where $\xi( s;\bbvarphi) =\theta( s;\bbvarphi_{0}) /\theta( s;\bbvarphi)
=\dsum_{j=0}^{\infty}\xi_{j}( \bbvarphi) s^{j}$. Then
%
%e2.19 #&#
\begin{equation} \label{at}
\inf_{\mathcal{T}_{3}}R_{n}( \bbtau) \geq
\mathop{\inf_{\lambda\in[ -\pi,\pi]}}_{\bbvarphi\in\Psi}%
\vert\xi( e^{i\lambda};\bbvarphi) \vert
^{2}\inf_{\delta\in\mathcal{I}_{3}}R_{n}( \bbtau^{\ast
}) -\sup_{\mathcal{T}_{3}}\frac{1}{n}\vert V_{n}( \bbtau) \vert.
\end{equation}
Assumption A1 implies [see (\ref{zf})]%
\[
\mathop{\inf_{\lambda\in[ -\pi,\pi]}}_{\bbvarphi%
\in\Psi}\vert\xi( e^{i\lambda};\bbvarphi)
\vert^{2}>\epsilon.
\]
Thus,
%
%e2.20 #&#
\begin{eqnarray} \label{cb}
\inf_{\mathcal{T}_{3}}R_{n}( \bbtau) &\geq&\epsilon
\inf_{\mathcal{I}_{3}}\frac{1}{n}\sum_{t=1}^{n}\Biggl(
\sum_{j=0}^{t-1}a_{j}\varepsilon_{t-j}\Biggr) ^{2}
\nonumber\\[-10pt]\\[-10pt]
&&{}-\sup_{\mathcal{T}_{3}}\frac{1}{n}\vert V_{n}( \bbtau%
) \vert-\sup_{\mathcal{I}_{3}}\frac{1}{n}\vert W_{n}(
\delta) \vert,\nonumber\vspace*{-2pt}
\end{eqnarray}
where $a_{j}=a_{j}( \delta_{0}-\delta) $, and by Lemma \ref{lemma2}
\[
W_{n}( \delta) =\epsilon\sum_{t=1}^{n}v_{t}^{2}(
\delta) +2\epsilon\sum_{t=1}^{n}v_{t}( \delta)
\sum_{j=0}^{t-1}a_{j}\varepsilon_{t-j}.\vspace*{-2pt}
\]
By Lemma \ref{lemma2} and (0.6) in the proof of Lemma \ref{lemma3}
in the supplementary material
\cite{hualde} (taking $\kappa=1/2$
there in both cases)%
%
%e2.21 #&#
\begin{equation} \label{cb2}
\sup_{\mathcal{I}_{3}}\frac{1}{n}\vert W_{n}( \delta) \vert=O_{p}\biggl(
n^{-1}+\frac{\log n}{n^{1/2}}\biggr) =o_{p}( 1) ,\vspace*{-2pt}
\end{equation}
and also by Lemma \ref{lemma3} (with $\kappa=1/2$ there)
%
%e2.22 #&#
\begin{equation} \label{cb1}
\sup_{\mathcal{T}_{3}}\frac{1}{n}\vert V_{n}( \bbtau%
) \vert=O_{p}\biggl( \frac{\log^{2}n}{n}\biggr) =o_{p}(
1) .\vspace*{-2pt}
\end{equation}
Next, note that for $\delta\in\mathcal{I}_{3}$%
%
%e2.23 #&#
\begin{equation} \label{51}
\frac{\partial a_{j}^{2}}{\partial\delta}=-2\bigl( \psi( j+\delta
_{0}-\delta) -\psi( \delta_{0}-\delta) \bigr)
a_{j}^{2}<0,\vspace*{-2pt}
\end{equation}
where we introduce the digamma function $\psi( x) =( d/dx)$log$\Gamma(
x)$.
From (\ref{51}) and the fact that $\psi( x) $ is strictly
increasing in $x>0$,
%
%e2.24 #&#
\begin{eqnarray} \label{h1}
\inf_{\mathcal{I}_{3}}n^{-1}\sum_{t=1}^{n}\Biggl(
\sum_{j=0}^{t-1}a_{j}\varepsilon_{t-j}\Biggr) ^{2}
&\geq& n^{-1}\sum_{t=1}^{n}\sum_{j=0}^{t-1}a_{j}^{2}\biggl(
\frac{1}{2}-\eta\biggr) \varepsilon_{t-j}^{2} \nonumber\\[-10pt]\\[-10pt]
&&{}-\sup_{\mathcal{I}_{3}}\Biggl\vert\frac{1}{n}\sum_{t=1}^{n}
\mathop{\sum\sum}^{t-1}_{j\neq k}a_{j}a_{k}\varepsilon
_{t-j}\varepsilon
_{t-k}\Biggr\vert.\nonumber\vspace*{-2pt}
\end{eqnarray}
By a very similar analysis to that of $( b) $ in (\ref{ai}), the
second term on the right-hand side of (\ref{h1}) is bounded by%
\begin{eqnarray*}
&&\frac{2}{n}\sup_{\mathcal{I}_{3}}\Biggl\vert
\dsum_{j=0}^{n-2}\dsum_{k=j+1}^{n-1}a_{j}a_{k}\dsum%
_{l=k-j+1}^{n-j}\varepsilon_{l}\varepsilon_{l-(k-j)}\Biggr\vert
\\[-3pt]
&&\qquad\leq\frac{2}{n}\sup_{\mathcal{I}_{3}}\Biggl\vert
\dsum_{j=0}^{n-2}a_{j}\dsum_{k=j+1}^{n-1}\dsum%
_{l=k-j+1}^{n-j}\varepsilon_{l}\varepsilon_{l-(k-j)}\Biggr\vert \\[-3pt]
&&\qquad\quad{}+\frac{2}{n}\sup_{\mathcal{I}_{3}}\Biggl\vert
\dsum_{j=0}^{n-2}a_{j}\dsum_{k=j+1}^{n-2}(a_{k+1}-a_{k})\dsum%
_{r=j+1}^{k}\dsum_{l=r-j+1}^{n-j}\varepsilon_{l}\varepsilon
_{l-(k-j)}\Biggr\vert,\vspace*{-2pt}\vadjust{\goodbreak}
\end{eqnarray*}
which has expectation bounded by
%
%e2.25 #&#
\begin{eqnarray} \label{h5}
&&\frac{K}{n^{{1/2}}}\sum_{j=1}^{n}j^{-{1/2}}+\frac{K}{n^{%
{1/2}}}\sum_{j=1}^{n}j^{-{1/2}}\sum_{k=1}^{n}%
( k+j) ^{-{3/2}}k^{{1/2}} \nonumber\\[-8pt]\\[-8pt]
&&\qquad\leq K\Biggl( 1+\frac{1}{n^{{1/2}}}\sum_{j=1}^{n}j^{-{1/2%
}-a}\sum_{k=1}^{n}k^{-1+a}\Biggr) \leq K\nonumber
\end{eqnarray}
for any $0<a<1/2$. Therefore, there exists a large enough $K$ such that
%
%e2.26 #&#
\begin{equation}\label{cc}
\Pr\Biggl( \sup_{\mathcal{I}_{3}}\Biggl\vert n^{-1}\sum_{t=1}^{n}%
\mathop{\sum\sum}^{t-1}_{j\neq k}a_{j}a_{k}\varepsilon
_{t-j}\varepsilon_{t-k}\Biggr\vert>K\Biggr) \rightarrow0
\end{equation}
as $n\rightarrow\infty$. Then, noting (\ref{cb}), (\ref{cb2}), (\ref
{cb1}%
), (\ref{cc}), we deduce (\ref{ac}) for $i=3$ if%
%
%e2.27 #&#
\begin{equation} \label{32}
\Pr\Biggl( \frac{1}{n}\sum_{t=1}^{n}\sum_{j=0}^{t-1}a_{j}^{2}%
\biggl( \frac{1}{2}-\eta\biggr) \varepsilon_{t-j}^{2}\leq K\Biggr)
\rightarrow0\qquad\mbox{as }n\rightarrow\infty.
\end{equation}
Now
\begin{eqnarray*}
\frac{1}{n}\sum_{t=1}^{n}\sum_{j=0}^{t-1}a_{j}^{2}\biggl(
\frac{1}{2}-\eta\biggr) \varepsilon_{t-j}^{2} &=&\sigma_{0}^{2}\frac{%
\Gamma( 2\eta) }{\Gamma^{2}(
{1/2}+\eta) }+\frac{1}{n}\sum_{t=1}^{n}\sum%
_{j=0}^{t-1}a_{j}^{2}\biggl( \frac{1}{2}-\eta \biggr) (
\varepsilon_{t-j}^{2}-\sigma_{0}^{2}) \\
&&{}-\frac{\sigma_{0}^{2}}{n}\sum_{t=1}^{n}\sum_{j=t}^{\infty
}a_{j}^{2}\biggl( \frac{1}{2}-\eta\biggr).
\end{eqnarray*}
The third term on the right is clearly $O( n^{-2\eta}) $, whereas, as
in the treatment of $( a) $ in (\ref{ai}%
), the second is $O_{p}( n^{-1/2}) $, so that (\ref{32}) holds as
$\Gamma( 2\eta) /\Gamma^{2}(
\frac{1}{2}+\eta) $ can be made arbitrarily large for small enough
$\eta$. This proves (\ref{ac}), and thus (\ref{new2}), for $i=3$.

\textit{Proof of} (\ref{new2}) \textit{for $i=2$}. Take $\eta<1/4$ and note
that $\mathcal{%
I}_{2}\subset\lbrack\delta_{0}-\kappa,\delta_{0}-\frac{1}{2}+\eta)$ for
$\kappa=\eta+ {\frac12} $. It follows from Lemma \ref{lemma2} and (0.6) in the
proof of Lemma \ref{lemma3} (see supplementary material \cite{hualde}) that
%
%e2.28 #&#
\begin{eqnarray} \label{ay}
\sup_{\mathcal{I}_{2}}\frac{1}{n}\vert W_{n}(\delta
)\vert&=&O_{p}\Biggl( \frac{1}{n}\dsum_{t=1}^{n}t^{2\eta-1}+%
\frac{1}{n}\dsum_{t=1}^{n}t^{\eta-%
{1/2}
}t^{\eta}\Biggr) \nonumber\\[-8pt]\\[-8pt]
&=&O_{p}( n^{2\eta-%
{1/2}}) =o_{p}(1).\nonumber
\end{eqnarray}
It follows from Lemma \ref{lemma3} that
%
%e2.29 #&#
\begin{equation} \label{az}
\sup_{\mathcal{T}_{2}}\frac{1}{n}\vert V_{n}(\bbtau%
)\vert=O_{p}( n^{2\eta-1}) =o_{p}(1).
\end{equation}
Denote\vspace*{1pt} $f_{n}( \delta) =n^{-1}\dsum_{t=1}^{n}(
\dsum_{j=0}^{t-1}a_{j}\varepsilon_{t-j}) ^{2}$. By
(\ref{ay}), (\ref{az}), it follows\break that~(\ref{ac}) for $i=2$ holds if for
arbitrarily large $K$%
%
%e2.30 #&#
\begin{equation} \label{ba}
\Pr\Bigl( \inf_{\mathcal{I}_{2}}f_{n}( \delta) >K\Bigr)
\rightarrow1
\end{equation}
as $n\rightarrow\infty$. Clearly,%
%
%e2.31 #&#
\begin{equation} \label{h2}
\inf_{\mathcal{I}_{2}}f_{n}( \delta) \geq\inf_{\mathcal{I}_{2}}%
\frac{n^{2( \delta_{0}-\delta) }}{n}\inf_{\mathcal{I}_{2}}\frac{%
1}{n^{2( \delta_{0}-\delta) }}\sum_{t=1}^{n}\Biggl(
\sum_{j=0}^{t-1}a_{j}\varepsilon_{t-j}\Biggr) ^{2}.
\end{equation}
Defining $b_{j,n}( d) =a_{j}( d) /n^{d-1}$, $%
b_{j,n}=b_{j,n}( \delta_{0}-\delta) $, the right-hand side of (\ref%
{h2}) is bounded below by%
%
%e2.32 #&#
\begin{equation} \label{h4}
\inf_{\mathcal{I}_{2}}\frac{1}{n^{2}}\sum_{j=0}^{n-1}b_{j,n}^{2}\sum%
_{l=1}^{n-j}\varepsilon_{l}^{2}-\sup_{\mathcal{I}_{2}}\frac{2}{n^{2}}%
\Biggl\vert
\sum_{j=0}^{n-2}\sum_{k=j+1}^{n-1}b_{j,n}b_{k,n}\sum%
_{l=k-j+1}^{n-j}\varepsilon_{l}\varepsilon_{l-( k-j)
}\Biggr\vert.\hspace*{-32pt}
\end{equation}
For $1\leq j\leq n$,
%
%e2.33 #&#
\begin{eqnarray} \label{h3}
\inf_{\mathcal{I}_{2}}b_{j,n}&\geq&\inf_{\mathcal{I}_{2}}\frac{\epsilon
}{%
\Gamma( \delta_{0}-\delta) }\inf_{\mathcal{I}_{2}}\biggl( \frac{%
j}{n}\biggr) ^{\delta_{0}-\delta-1}\geq\frac{\epsilon}{\Gamma(
{1/2}+\eta) }\biggl( \frac{j}{n}\biggr) ^{\eta-{1/2}},\hspace*{-30pt}
\nonumber\\
\sup_{\mathcal{I}_{2}}b_{j,n}&\leq&\sup_{\mathcal{I}_{2}}\frac{K}{\Gamma
( \delta_{0}-\delta) }\sup_{\mathcal{I}_{2}}\biggl( \frac{j}{n}%
\biggr) ^{\delta_{0}-\delta-1}\leq\frac{K}{\sqrt{\pi}}\biggl( \frac{j}{n}%
\biggr) ^{-{1/2}}.\hspace*{-30pt}
\end{eqnarray}
Then by (\ref{h3}), using summation by parts as in the analysis of $(
b) $ in (\ref{ai}), the expectation of the second term in (\ref{h4})
is bounded by%
\[
\frac{K}{n}\sum_{j=1}^{n}\biggl( \frac{j}{n}\biggr) ^{-{1/2}}+%
\frac{K}{n^{{1/2}}}\sum_{j=1}^{n}j^{-{1/2}%
}\sum_{k=1}^{n}k^{{1/2}}( k+j) ^{-{3/2}},
\]
which, noting (\ref{h5}), is $O( 1) $. Next, the first term in (%
\ref{h4}) is bounded below by%
%
%e2.34 #&#
\begin{equation} \label{h6}\qquad
\frac{\sigma_{0}^{2}}{n^{2}}\sum_{j=0}^{n-1}( n-j)
b_{j,n}^{2}( 1/2+\eta) -\frac{1}{n^{2}}\sum%
_{j=0}^{n-1}b_{j,n}^{2}( 1/2) \Biggl\vert
\sum_{l=1}^{n-j}( \varepsilon_{l}^{2}-\sigma_{0}^{2})
\Biggr\vert.
\end{equation}
Using (\ref{h3}) it can be easily shown that the second term in (\ref{h6})
is\break $O_{p}( n^{-{3/2}}\times\sum_{j=1}^{n}\frac{n}{j}) =O_{p}(n^{-%
{1/2}}\log n)$, whereas the first term is bounded below by
%
%e2.35 #&#
\begin{eqnarray} \label{h7}
&&\frac{\epsilon}{n}\dsum_{j=1}^{n}\biggl\{ \biggl( \frac{j}{n}\biggr)
^{2\eta-1}-\biggl( \frac{j}{n}\biggr) ^{2\eta}\biggr\} \nonumber\\
&&\qquad\geq\frac{\epsilon}{2}\dint_{1/n}^{1}\{ x^{2\eta
-1}-x^{2\eta}\} \,dx=\frac{\epsilon}{2}\biggl[ \frac{x^{2\eta}}{2\eta}%
-\frac{x^{2\eta+1}}{2\eta+1}\biggr] _{1/n}^{1} \\
&&\qquad=\frac{\epsilon}{4\eta(2\eta+1)}-O_{p}(n^{-2\eta}).\nonumber
\end{eqnarray}
Then (\ref{ba}) holds because the right-hand side of (\ref{h7}) can be made
arbitrarily large on setting $\eta$ arbitrarily close to zero. This
proves (\ref{ac}), and thus (\ref{new2}), for $i=2$.

\textit{Proof of} (\ref{new2}) \textit{for $i=1$}. Noting that $R_{n}( \bbtau%
) \geq n^{-2}( \dsum_{t=1}^{n}\varepsilon_{t}(
\bbtau) ) ^{2}$,
%
%e2.36 #&#
\begin{equation}\label{q1}
\Pr\Bigl( \inf_{\mathcal{T}_{1}}R_{n}( \bbtau) >K\Bigr)
\geq\Pr\Biggl( n^{2\eta}\inf_{\mathcal{T}_{1}}\Biggl( \frac{1}{n^{\delta
_{0}-\delta+{1/2}}}\dsum_{t=1}^{n}\varepsilon_{t}(
\bbtau) \Biggr) ^{2}>K\Biggr) ,\hspace*{-30pt}
\end{equation}
because $\delta_{0}-\delta\geq\frac{1}{2}+\eta$. Clearly $\dsum%
_{t=1}^{n}\varepsilon_{t}( \bbtau)
=\dsum_{j=0}^{n-1}d_{j}( \bbtau) u_{n-j}$, where
\[
d_{j}( \bbtau) =\dsum_{k=0}^{j}c_{k}( \bbtau) =\dsum_{k=0}^{j}\phi_{k}( \bbvarphi%
) \dsum_{l=0}^{j-k}a_{l}( \delta_{0}-\delta)
=\dsum_{k=0}^{j}\phi_{k}( \bbvarphi)
a_{j-k}( \delta_{0}-\delta+1) .
\]
For arbitrarily small $\epsilon>0$, the right-hand side of (\ref{q1})
is bounded
from below by
%
%e2.37 #&#
\begin{equation} \label{q2}
\Pr\Biggl( \inf_{\mathcal{T}_{1}}\Biggl( \frac{1}{n^{\delta_{0}-\delta+%
{1/2}}}\dsum_{t=1}^{n}\varepsilon_{t}( \bbtau%
) \Biggr) ^{2}>\epsilon\Biggr)
\end{equation}
for $n$ large enough, so it suffices to show (\ref{q2}) $\rightarrow1$
as $%
n\rightarrow\infty$. First%
\[
\frac{1}{n^{\delta_{0}-\delta+{1/2}}}\dsum_{t=1}^{n}%
\varepsilon_{t}( \bbtau) =\phi( 1;\bbvarphi%
) \theta( 1;\bbvarphi_{0}) h_{n}( \delta
) +r_{n}( \bbtau) ,
\]
where $h_{n}( \delta)
=n^{-1/2}\dsum_{j=0}^{n-1}b_{j,n}( \delta_{0}-\delta
+1) \varepsilon_{n-j}$, $b_{j,n}( \cdot) $ was defined
below~(\ref{h2}), and%
%
%e2.38 #&#
\begin{eqnarray} \label{ee}
r_{n}( \bbtau) &=&-\frac{1}{n^{{1/2}}}%
\dsum_{j=0}^{n-1}b_{j,n}( \delta_{0}-\delta+1)
\dsum_{k=j+1}^{\infty}\phi_{k}( \bbvarphi)
u_{n-j}\nonumber\\
&&{}-\frac{1}{n^{{1/2}}}\dsum_{j=1}^{n-1}s_{j,n}(
\bbtau) u_{n-j} \\
&&{}+\frac{\phi( 1;\bbvarphi) }{n^{{1/2}}}%
\dsum_{j=0}^{n-1}b_{j,n}( \delta_{0}-\delta+1) \bigl(
u_{n-j}-\theta( 1;\bbvarphi_{0}) \varepsilon
_{n-j}\bigr) \nonumber
\end{eqnarray}
for
\[
s_{j,n}( \bbtau) =\dsum_{k=0}^{j-1}\bigl(
b_{k+1,n}( \delta_{0}-\delta+1) -b_{k,n}( \delta
_{0}-\delta+1) \bigr) \dsum_{l=0}^{k}\phi_{j-l}(
\bbvarphi) ,
\]
where (\ref{ee}) is routinely derived, noting that by summation by parts
\begin{eqnarray*}
d_{j}( \bbtau) &=&a_{j}( \delta_{0}-\delta+1)\\
&&{}\times
\dsum_{k=0}^{j}\phi_{k}( \bbvarphi)
-\dsum_{k=0}^{j-1}\bigl( a_{k+1}( \delta_{0}-\delta+1)
-a_{k}( \delta_{0}-\delta+1) \bigr) \dsum_{l=0}^{k}\phi
_{j-l}( \bbvarphi) .
\end{eqnarray*}
Now%
\begin{eqnarray*}
\inf_{\mathcal{T}_{1}}\Biggl( \frac{1}{n^{\delta_{0}-\delta+{1/2}}}%
\dsum_{t=1}^{n}\varepsilon_{t}( \bbtau) \Biggr)
^{2} &\geq&\theta^{2}( 1;\bbvarphi_{0}) \inf_{\bbPsi}\phi^{2}( 1;\bbvarphi) \inf_{\mathcal{I}%
_{1}}h_{n}^{2}( \delta) \\
&&{}-K\sup_{\bbPsi}\vert\phi( 1;\bbvarphi)
\vert\sup_{\mathcal{I}_{1}}\vert h_{n}( \delta)
\vert\sup_{\mathcal{T}_{1}}\vert r_{n}( \bbtau%
) \vert.
\end{eqnarray*}
Noting (\ref{zf}) and that under A1, $\sup_{\bbPsi}\vert\phi
( 1;\bbvarphi) \vert<\infty$, the required result
follows on showing that%
%
%e2.39 #&#
%e2.40 #&#
%e2.41 #&#
\begin{eqnarray}
\label{qqq1}
\sup_{\mathcal{T}_{1}}\vert r_{n}( \bbtau)
\vert&=&o_{p}( 1) ,\\
\label{qqq2}
\sup_{\mathcal{I}_{1}}\vert h_{n}( \delta) \vert
&=&O_{p}( 1) ,\\
\label{qqq3}
\Pr\Bigl( \inf_{\mathcal{I}_{1}}h_{n}^{2}( \delta) >\epsilon
\Bigr) &\rightarrow&1
\end{eqnarray}
as $n\rightarrow\infty$.

The proof of (\ref{qqq2}) is omitted as it is similar to and much easier
than the proof of (\ref{qqq1}), which we now give. Let $r_{n}( \bbtau) =\dsum_{i=1}^{3}r_{in}( \bbtau) $.
By the Cauchy inequality%
\[
\sup_{\mathcal{T}_{1}}\vert r_{1n}( \bbtau)
\vert\leq\frac{1}{n^{{1/2}}}\Biggl(
\dsum_{j=0}^{n-1}\sup_{\mathcal{I}_{1}}b_{j,n}^{2}( \delta
_{0}-\delta+1) \Biggl( \sup_{\Psi}\dsum_{k=j+1}^{\infty
}\vert\phi_{k}( \bbvarphi) \vert\Biggr)
^{2}\dsum_{j=1}^{n}u_{j}^{2}\Biggr) ^{{1/2}},
\]
so that by (\ref{eee}), noting that $E(
\dsum_{j=1}^{n}u_{j}^{2}) ^{1/2}\leq Kn^{1/2}$,
\begin{eqnarray*}
{E\sup_{\mathcal{T}_{1}}}\vert r_{1n}( \bbtau)
\vert &\leq&K\Biggl( \dsum_{j=1}^{n}\sup_{\mathcal{I}%
_{1}}\biggl( \frac{j}{n}\biggr) ^{2( \delta_{0}-\delta) }\Biggl(
\dsum_{k=j+1}^{\infty}k^{-1-\varsigma}\Biggr) ^{2}\Biggr) ^{{1%
/2}} \\
&\leq&K\Biggl( \dsum_{j=1}^{n}\biggl( \frac{j}{n}\biggr) ^{1+2\eta
}j^{-2\varsigma}\Biggr) ^{{1/2}}\leq Kn^{{1/2}-\varsigma
}=o( 1) ,
\end{eqnarray*}
because $\varsigma>1/2$ by A1(ii). Next, by summation by parts%
\[
r_{2n}( \bbtau) =-\frac{s_{n-1,n}( \bbtau%
) }{n^{{1/2}}}\dsum_{j=1}^{n-1}u_{n-j}+\frac{1}{n^{{1%
/2}}}\dsum_{j=1}^{n-2}\bigl( s_{j+1,n}( \bbtau)
-s_{j,n}( \bbtau) \bigr) \dsum_{k=1}^{j}u_{n-k},
\]
so%
%
%e2.42 #&#
\begin{eqnarray} \label{q11}
\sup_{\mathcal{T}_{1}}\vert r_{2n}( \bbtau)
\vert &\leq&\frac{{\sup_{\mathcal{T}_{1}}}\vert s_{n-1,n}(
\bbtau) \vert}{n^{{1/2}}}\Biggl\vert
\dsum_{j=1}^{n-1}u_{n-j}\Biggr\vert \nonumber\\[-8pt]\\[-8pt]
&&{}+\frac{1}{n^{{1/2}}}\dsum_{j=1}^{n-2}\sup_{\mathcal{T}%
_{1}}\vert s_{j+1,n}( \bbtau) -s_{j,n}( \bbtau) \vert\Biggl\vert
\dsum_{k=1}^{j}u_{n-k}\Biggr\vert.\nonumber
\end{eqnarray}
Given that $a_{k+1}( \delta_{0}-\delta+1) -a_{k}( \delta
_{0}-\delta+1) =a_{k+1}( \delta_{0}-\delta)$,
\[
s_{j,n}( \bbtau) =\frac{1}{n^{\delta_{0}-\delta}}%
\dsum_{k=0}^{j-1}a_{k+1}( \delta_{0}-\delta)
\dsum_{l=0}^{k}\phi_{j-l}( \bbvarphi) ,
\]
so as $E\vert\dsum_{j=1}^{n-1}u_{j}\vert\leq Kn^{1/2}$%
, noting (\ref{eee}) and Stirling's approximation, the expectation of the
first term on the right-hand side of (\ref{q11}) is bounded by%
\begin{eqnarray*}
&&K\dsum_{k=1}^{n}\sup_{\mathcal{I}_{1}}\biggl( \frac{k}{n}\biggr)
^{\delta_{0}-\delta}k^{-1}\dsum_{l=1}^{k}( n-l)
^{-1-\varsigma} \\
&&\qquad\leq\frac{K}{n^{{1/2}+\eta}}\dsum%
_{k=1}^{n}k^{-{1/2}+\eta}( n-k) ^{-{1/2}} \\
&&\qquad\leq\frac{K}{n^{{1/2}}}\frac{1}{n}\dsum_{k=1}^{n}\biggl(
\frac{k}{n}\biggr) ^{-{1/2}+\eta}\biggl( 1-\frac{k}{n}\biggr) ^{-{%
1/2}}\\
&&\qquad\leq Kn^{-{1/2}}.
\end{eqnarray*}
Next, noting that $a_{j+1}( \delta_{0}-\delta) -a_{j}(
\delta_{0}-\delta) =a_{j+1}( \delta_{0}-\delta-1) $, it
can be shown that%
%
%e2.43 #&#
\begin{eqnarray}\label{q12}
s_{j+1,n}( \bbtau) -s_{j,n}( \bbtau)
&=&\frac{1}{n^{\delta_{0}-\delta}}\dsum_{k=1}^{j}\phi_{k}(
\bbvarphi) \dsum_{l=j-k+2}^{j+1}a_{l}( \delta
_{0}-\delta-1) \nonumber\\[-8pt]\\[-8pt]
&&{}+\frac{\phi_{j+1}( \bbvarphi) }{n^{\delta_{0}-\delta}%
}\dsum_{l=1}^{j+1}a_{l}( \delta_{0}-\delta) .
\nonumber
\end{eqnarray}
Thus, noting that, uniformly in $j$, $n$, $E\vert
{\dsum_{k=1}^{j}u_{n-k}}\vert\leq Kj^{1/2}$, by previous
arguments the contribution of the last term on the right-hand side of
(\ref{q12})
to the expectation of the second term on the right-hand side of (\ref
{q11}) is
bounded by%
\begin{eqnarray*}
\frac{K}{n^{{1/2}}}\dsum_{j=1}^{n}j^{{1/2}%
}j^{-1-\varsigma}\sup_{\mathcal{I}_{1}}\biggl( \frac{j}{n}\biggr) ^{\delta
_{0}-\delta}&\leq&\frac{K}{n^{{1/2}}}\dsum_{j=1}^{n}j^{-{1%
/2}-\varsigma}\biggl( \frac{j}{n}\biggr) ^{{1/2}+\eta}\\&\leq&
Kn^{-\varsigma}.
\end{eqnarray*}
By identical arguments, the contribution of the first term on the
right-hand side
of (\ref{q12}) to the expectation of the last term on the right-hand
side of~(\ref{q11}) is bounded by
%
%e2.44 #&#
\begin{eqnarray} \label{aa}
&&\frac{K}{n^{{1/2}}}\dsum_{j=1}^{n}j^{{1/2}%
}\dsum_{k=1}^{j-1}k^{-1-\varsigma}\dsum_{l=j-k}^{j}\sup_{%
\mathcal{I}_{1}}\biggl( \frac{l}{n}\biggr) ^{\delta_{0}-\delta}l^{-2}
\nonumber\\[-8pt]\\[-8pt]
&&\qquad\leq\frac{K}{n^{1+\eta}}\dsum_{j=1}^{n}j^{{1/2}%
}\dsum_{k=1}^{j-1}k^{-1-\varsigma}\dsum_{l=j-k}^{j}l^{-{3%
/2}+\eta}.\nonumber
\end{eqnarray}
Given that $\dsum_{l=j-k}^{j}l^{-{3/2}+\eta}\leq K(
j-k) ^{-{3/2}+\eta}k$, the right-hand side of (\ref{aa}) is bounded
by%
%
%e2.45 #&#
\begin{eqnarray} \label{as}
&&\frac{K}{n^{1+\eta}}\dsum_{j=1}^{n}j^{{1/2}%
}\dsum_{k=1}^{j-1}k^{-\varsigma}( j-k) ^{-{3/2}%
+\eta} \nonumber\\
&&\qquad\leq\frac{K}{n^{1+\eta}}\dsum_{j=1}^{n}j^{{1/2}%
}\dsum_{k=1}^{[ j/2] }k^{-\varsigma}( j-k) ^{-%
{3/2}+\eta} \\
&&\qquad\quad{}+\frac{K}{n^{1+\eta}}\dsum_{j=1}^{n}j^{{1/2}}\dsum_{k=%
[ j/2] +1}^{j-1}k^{-\varsigma}( j-k) ^{-{3/2}%
+\eta},\nonumber
\end{eqnarray}
where $[ \cdot] $ denotes integer part. Clearly, the right-hand side
of (\ref{as}) is bounded by%
\[
\frac{K}{n^{1+\eta}}\dsum_{j=1}^{n}j^{{1/2}}\Biggl( j^{-{3%
/2}+\eta}j^{1-\varsigma}+j^{-\varsigma}\dsum_{k=1}^{\infty}k^{-%
{3/2}+\eta}\Biggr) \leq K( n^{-\varsigma}+n^{{1/2}%
-\varsigma-\eta}) ,
\]
so $\sup_{\mathcal{T}_{1}}\vert r_{2n}( \bbtau)
\vert=o_{p}( 1) $ because $\varsigma>1/2$. Next, writing $%
u_{t}=\theta( 1;\bbvarphi_{0}) \varepsilon_{t}+
\widetilde{\varepsilon}_{t-1}-\widetilde{\varepsilon}_{t}$, for $%
\widetilde{\varepsilon}_{t}=\dsum_{j=0}^{\infty}\widetilde{\theta
}_{j}( \bbvarphi_{0}) \varepsilon_{t-j}$, $\widetilde{%
\theta}_{j}( \bbvarphi_{0})
=\dsum_{k=j+1}^{\infty}\theta_{k}( \bbvarphi%
_{0}) $, where, by A1, A2, $\widetilde{\varepsilon}_{t}$ is well
defined in the mean square sense, we have%
\[
r_{3n}( \bbtau) =-\frac{\phi( 1;\bbvarphi%
) }{n^{\delta_{0}-\delta+{1/2}}}\Biggl(
\dsum_{j=0}^{n-1}a_{j}( \delta_{0}-\delta) \widetilde{%
\varepsilon}_{n-k}-a_{n-1}( \delta_{0}-\delta+1) \widetilde{%
\varepsilon}_{0}\Biggr) .
\]
In view of previous arguments, it is straightforward to show that\break ${\sup
_{%
\mathcal{T}_{1}}}\vert r_{3n}( \bbtau) \vert
=o_{p}( 1) $, to conclude the proof of (\ref{qqq1}).

Finally, we prove (\ref{qqq3}). Considering $h_{n}( \delta) $ as
a process indexed by $\delta$, we show first that%
%
%e2.46 #&#
\begin{equation} \label{q14}
h_{n}( \delta) \Rightarrow\dint_{0}^{1}\frac{(
1-s) ^{\delta_{0}-\delta}}{\Gamma( \delta_{0}-\delta
+1) }\,dB( s) ,
\end{equation}
where $B( s) $ is a scalar Brownian motion with variance $\sigma
_{0}^{2}$ and $\Rightarrow$ means weak convergence in the space of
continuous functions on $\mathcal{I}_{1}$. We give this space the uniform
topology. Convergence of the finite-dimensional distributions follows by
Theorem 1 of \cite{hosoya2}, noting that A2 implies conditions A(i), A(ii)
and A(iii) in \cite{hosoya2} (in particular A2 implies that the fourth-order
cumulant spectral density function of $\varepsilon_{t}$ is bounded). Next,
by Theorem 12.3 of \cite{billingsley}, if for all fixed $\delta\in
\mathcal{%
I}_{1}$ $h_{n}( \delta) $ is a tight sequence, and if for all $%
\delta_{1},\delta_{2}\in\mathcal{I}_{1}$ and for~$K$ not depending
on $%
\delta_{1},\delta_{2},n$%
%
%e2.47 #&#
\begin{equation} \label{q6}
E\bigl( h_{n}( \delta_{1}) -h_{n}( \delta_{2})
\bigr) ^{2}\leq K( \delta_{1}-\delta_{2}) ^{2},
\end{equation}
then the process $h_{n}( \delta) $ is tight, and (\ref{q14})
would follow. First, for fixed $\delta$, it is straightforward to show that
$\sup_{n}E( h_{n}^{2}( \delta) ) <\infty$, so $%
h_{n}( \delta) $ is uniformly integrable and therefore tight.
Next,%
\begin{eqnarray*}
&&E\bigl( h_{n}( \delta_{1}) -h_{n}( \delta_{2})
\bigr) ^{2}\\
&&\qquad=\frac{\sigma_{0}^{2}}{n}\dsum_{j=0}^{n-1}\bigl(
b_{j,n}( \delta_{0}-\delta_{1}+1) -b_{j,n}( \delta
_{0}-\delta_{2}+1) \bigr) ^{2} \\
&&\qquad=\frac{\sigma_{0}^{2}( \delta_{1}-\delta_{2}) ^{2}}{n}%
\dsum_{j=0}^{n-1}\frac{( a_{j}^{\prime}( \delta_{0}-%
\overline{\delta}+1) -a_{j}( \delta_{0}-\overline{\delta}%
+1) \log n) ^{2}}{n^{2( \delta_{0}-\overline{\delta}%
) }}
\end{eqnarray*}
by the mean value theorem, where $\overline{\delta}=\overline{\delta}_{n}$
is an intermediate point between~$\delta_{1}$ and~$\delta_{2}$. As in
Lemma D.1 of \cite{robhualde},%
\begin{eqnarray*}
&&a_{j}^{\prime}( \delta_{0}-\overline{\delta}+1) -a_{j}(
\delta_{0}-\overline{\delta}+1) \log n \\
&&\qquad=\bigl( \psi( j+\delta_{0}-\overline{\delta}+1) -\psi(
\delta_{0}-\overline{\delta}+1) -\log n\bigr) a_{j}( \delta
_{0}-\overline{\delta}+1) .
\end{eqnarray*}
Now (\ref{q6}) holds on showing that, for $\overline{\delta}\in
\mathcal{I}%
_{1}$,%
%
%e2.48 #&#
%e2.49 #&#
\begin{eqnarray}
 \label{q8}
\frac{\psi^{2}( \delta_{0}-\overline{\delta}+1) }{n}%
\dsum_{j=0}^{n-1}b_{j,n}^{2}( \delta_{0}-\overline{\delta}%
+1) &\leq&K,\\
\label{q9}
\frac{1}{n}\dsum_{j=0}^{n-1}\bigl( \psi( j+\delta_{0}-%
\overline{\delta}+1) -\log n\bigr) ^{2}b_{j,n}^{2}( \delta_{0}-%
\overline{\delta}+1) &\leq&K.
\end{eqnarray}
By Stirling's approximation, the left-hand side of (\ref{q8}) is
bounded by%
\begin{eqnarray*}
&&
K\frac{\psi^{2}( \delta_{0}-\bigtriangledown_{1}+1) }{n}%
\dsum_{j=1}^{n}\sup_{\mathcal{I}_{1}}\biggl( \frac{j}{n}\biggr)
^{2( \delta_{0}-\delta) }\\
&&\qquad\leq K\frac{\psi^{2}( \delta
_{0}-\bigtriangledown_{1}+1) }{n}\dsum_{j=1}^{n}\sup_{\mathcal{%
I}_{1}}\biggl( \frac{j}{n}\biggr) ^{1+2\eta}\leq K.
\end{eqnarray*}
Regarding (\ref{q9}), it can be shown that
uniformly in $\mathcal{I}_{1}$, $\psi( j+\delta_{0}-\overline{\delta}
+1) =\log j+O( j^{-1}) $ (see, e.g., \cite{abramowitz},
page 259). Thus, apart from a remainder term of smaller order, the
left-hand side of
(\ref{q9}) is bounded by
%
%e2.50 #&#
\begin{equation} \label{q10}
K\frac{1}{n}\dsum_{j=1}^{n}\biggl( \log\frac{j}{n}\biggr)
^{2}b_{j,n}^{2}( \delta_{0}-\overline{\delta}+1) \leq K\frac{1}{%
n}\dsum_{j=1}^{n}\biggl( \log\frac{j}{n}\biggr) ^{2}\biggl( \frac{j}{n}%
\biggr) ^{1+2\eta}\hspace*{-30pt}
\end{equation}
uniformly in $\mathcal{I}_{1}$, the right-hand side of (\ref{q10})
being bounded
by\break $K\dint_{0}^{1}( \log x) ^{2}\,dx=2K$, to conclude the
proof of tightness. Then by the continuous mapping theorem%
\[
\inf_{\mathcal{I}_{1}}h_{n}^{2}( \delta) \rightarrow_{d}\inf_{%
\mathcal{I}_{1}}\biggl( \dint_{0}^{1}\frac{( 1-s) ^{\delta
_{0}-\delta}}{\Gamma( \delta_{0}-\delta+1) }\,dB( s)
\biggr) ^{2}.
\]
This is a.s. positive because the quantity whose infimum is taken is a
$\chi
_{1}^{2}$ random variable times $\sigma_{0}^{2}/[ \{ 2(
\delta_{0}-\delta) +1\} \Gamma( \delta_{0}-\delta
+1) ^{2}] $, which is bounded away from zero on $\mathcal{I}_{1}$.
Thus as $n\rightarrow\infty$
\[
\Pr\Bigl( \inf_{\mathcal{I}_{1}}h_{n}^{2}( \delta) >\epsilon
\Bigr) \rightarrow\Pr\biggl( \inf_{\mathcal{I}_{1}}\biggl(
\dint_{0}^{1}\frac{( 1-s) ^{\delta_{0}-\delta}}{\Gamma
( \delta_{0}-\delta+1) }\,dB( s) \biggr) ^{2}>\epsilon
\biggr) ,
\]
and (\ref{qqq3}) follows as $\epsilon$ is arbitrarily small. Then we
conclude (\ref{ac}), and thus (\ref{new2}), for $i=1$.
\end{pf}

%s2.2 #&#
\subsection{\texorpdfstring{Asymptotic normality of $\widehat{\bbtau}$}{Asymptotic normality of tau}}

This requires an additional regularity condition.

\begin{enumerate}[A3.]
\item[A3.]

\begin{enumerate}[(iii)]
\item[(i)]
\[
\bbtau_{0}\in \operatorname{int}\mathcal{T};
\]

\item[(ii)] for all $\lambda$, $\theta( e^{i\lambda};\bbvarphi) $ is twice continuously
differentiable in $\bbvarphi$ on
a~closed neighborhood $\mathcal{N}_{\epsilon}( \bbvarphi%
_{0}) $ of radius $0<\epsilon<1/2$ about $\bbvarphi_{0};$

\item[(iii)] the matrix%
\[
\mathbf{A}=\pmatrix{
\pi^{2}/6 & \displaystyle -\sum_{j=1}^{\infty}\mathbf{b}_{j}^{\prime}(
\bbvarphi_{0}) /j \vspace*{2pt}\cr
\displaystyle -\sum_{j=1}^{\infty}\mathbf{b}_{j}( \bbvarphi%
_{0}) /j & \displaystyle \sum_{j=1}^{\infty}\mathbf{b}_{j}( \bbvarphi_{0}) \mathbf{b}_{j}^{\prime}( \bbvarphi%
_{0})}
\]
is nonsingular, where $\mathbf{b}_{j}( \bbvarphi_{0})
=\dsum_{k=0}^{j-1}\theta_{k}( \bbvarphi_{0})
\partial\phi_{j-k}( \bbvarphi_{0}) /\partial\bbvarphi$.
\end{enumerate}
\end{enumerate}

By compactness of $\mathcal{N}_{\epsilon}( \bbvarphi%
_{0}) $ and continuity of $\partial\phi_{j}( \bbvarphi%
) /\partial\varphi_{i}$, $\partial^{2}\phi_{j}( \bbvarphi) /\partial\varphi_{i}\,
\partial\varphi_{l}$, for all $j$,
with $i,l=1,\ldots,p$, where $\varphi_{i}$ is the $i$th element of
$\bbvarphi$, A1(ii), A1(iv) and A3(ii) imply that, as $j\rightarrow
\infty$
\[
\sup_{\bbvarphi\in\mathcal{N}_{\epsilon}( \bbvarphi%
_{0}) }\biggl\vert\frac{\partial\phi_{j}( \bbvarphi%
) }{\partial\varphi_{i}}\biggr\vert=O\bigl( j^{-( 1+\varsigma
) }\bigr), \qquad\sup_{\bbvarphi\in
\mathcal{N}_{\epsilon}( \bbvarphi_{0}) }\biggl\vert\frac{%
\partial^{2}\phi_{j}( \bbvarphi) }{\partial\varphi
_{i}\,\partial\varphi_{l}}\biggr\vert=O\bigl( j^{-( 1+\varsigma
) }\bigr),
\]
which again is satisfied in the ARMA case. As with A1, A3 is similar to
conditions employed under stationarity, and can readily be checked in
general.

%th2.2 #&#
\begin{theorem}\label{theo2.2}
Let (\ref{a}), (\ref{b}) and \textup{A1--A3} hold. Then as $n\rightarrow
\infty$%
%
%e2.51 #&#
\begin{equation} \label{213}
n^{{1/2}}( \widehat{\bbtau}-\bbtau_{0})
\rightarrow_{d}N( 0,\mathbf{A}^{-1}) .\vadjust{\goodbreak}
\end{equation}
\end{theorem}
\begin{pf}
The proof standardly involves use of the mean value theorem,
approximation of a score function by a martingale
so as to apply a martingale convergence theorem, and convergence in
probability of a
Hessian in a neighborhood of $\bbtau_{0}$. From the mean value
theorem, (\ref{213}) follows if
we prove that
%
%e2.52 #&#
\begin{eqnarray} \label{x2}
\frac{\sqrt{n}}{2}\frac{\partial R_{n}( \bbtau_{0}) }{%
\partial\bbtau}&\rightarrow_{d}& N( 0,\sigma_{0}^{4}%
\mathbf{A}) ,
\\
%
%
%e2.53 #&#
\label{x3}
\frac{1}{2}\frac{\partial^{2}R_{n}( \overline{\bbtau}) }{%
\partial\bbtau\partial\bbtau^{\prime}}&\rightarrow
_{p}&
\sigma_{0}^{2}\mathbf{A},
\end{eqnarray}
where $\Vert\overline{\bbtau}-\bbtau_{0}\Vert
\leq\Vert\widehat{\bbtau}-\bbtau_{0}\Vert$.

\textit{Proof of} (\ref{x2}). It suffices to prove
%
%e2.54 #&#
\begin{equation}\label{x4}
\frac{\sqrt{n}}{2}\frac{\partial R_{n}( \bbtau_{0}) }{%
\partial\bbtau}-\frac{1}{\sqrt{n}}\sum_{t=2}^{n}\varepsilon
_{t}\sum_{j=1}^{\infty}\mathbf{m}_{j}( \bbvarphi%
_{0}) \varepsilon_{t-j}=o_{p}( 1)
\end{equation}
and
%
%e2.55 #&#
\begin{equation} \label{x5}
\frac{1}{\sqrt{n}}\sum_{t=2}^{n}\varepsilon
_{t}\sum_{j=1}^{\infty}\mathbf{m}_{j}( \bbvarphi%
_{0}) \varepsilon_{t-j}\rightarrow_{d}N( 0,\sigma_{0}^{4}%
\mathbf{A}) ,
\end{equation}
where $\mathbf{m}_{j}( \bbvarphi_{0}) =( -j^{-1},%
\mathbf{b}_{j}^{\prime}( \bbvarphi_{0}) ) ^{\prime
}$. By Lemma \ref{lemma2}, the left-hand side of (\ref{x4}) is the $(p+1)\times1$
vector $( r_{1}+r_{2}+r_{3},( \mathbf{s}_{1}+\mathbf{s}_{2})
^{\prime}) ^{\prime}$, where
\begin{eqnarray*}
r_{1} &=&\frac{1}{\sqrt{n}}\sum_{t=2}^{n}\varepsilon
_{t}\sum_{j=t}^{\infty}\frac{1}{j}\varepsilon_{t-j},\\
r_{2}&=&\frac{1}{\sqrt{n}}\sum_{t=2}^{n}\varepsilon
_{t}\sum_{j=1}^{t-1}\frac{1}{j}\sum_{k=t-j}^{\infty}\phi
_{k}( \bbvarphi_{0}) u_{t-j-k}, \\
r_{3} &=&-\frac{1}{\sqrt{n}}\sum_{t=2}^{n}v_{t}( \delta
_{0}) \sum_{j=1}^{t-1}\frac{1}{j}\sum_{k=0}^{t-j-1}\phi
_{k}( \bbvarphi_{0}) u_{t-j-k}, \\
\mathbf{s}_{1}&=&\frac{1}{\sqrt{n}}\sum_{t=2}^{n}\varepsilon
_{t}\sum_{j=t}^{\infty}\frac{\partial\phi_{j}( \bbvarphi_{0}) }{\partial\bbvarphi}u_{t-j}, \\
\mathbf{s}_{2} &=&\frac{1}{\sqrt{n}}\sum_{t=2}^{n}v_{t}( \delta
_{0}) \sum_{j=1}^{t-1}\frac{\partial\phi_{j}( \bbvarphi_{0}) }{\partial\bbvarphi}u_{t-j}.
\end{eqnarray*}
Clearly, $E( r_{1}) =0$, and
\[
\operatorname{Var}( r_{1}) =\frac{1}{n}\sum_{t=2}^{n}\sum_{j=t}^{%
\infty}\sum_{s=2}^{n}\sum_{k=s}^{\infty}\frac{1}{jk}E(
\varepsilon_{t}\varepsilon_{s}\varepsilon_{t-j}\varepsilon_{s-k}) =%
\frac{\sigma_{0}^{4}}{n}\sum_{t=2}^{n}\sum_{j=t}^{\infty}%
\frac{1}{j^{2}}=O\biggl( \frac{\log n}{n}\biggr) ,
\]
noting that, by A2, the $\varepsilon_{t}$ and $\varepsilon
_{t}^{2}-\sigma
_{0}^{2}$ are martingale difference sequences. Thus, $r_{1}=O_{p}(
n^{-1/2}\log^{1/2}n) $. Next, $E( r_{2}) =0$, and $%
\operatorname{Var}( r_{2}) $ equals%
%
%e2.56 #&#
\begin{equation} \label{x6}
\frac{1}{n}\sum_{t=2}^{n}\sum_{j=1}^{t-1}\sum_{k=t-j}^{%
\infty
}\sum_{s=2}^{n}\sum_{l=1}^{s-1}\sum_{m=s-l}^{\infty}%
\frac{\phi_{k}( \bbvarphi_{0}) \phi_{m}
( \bbvarphi_{0}) }{jl}E( \varepsilon_{t}\varepsilon
_{s}u_{t-j-k}u_{s-l-m}) .\hspace*{-35pt}
\end{equation}
From (\ref{b}) and A2, the expectation is $\sigma_{0}^{2}\gamma
_{j+k-\ell
-m}$ for $s=t$, and zero otherwise. By A1, $u_{t}$ has bounded spectral
density. Thus, (\ref{x6}) is bounded by
\begin{eqnarray*}
&&
K\frac{1}{n}\sum_{t=2}^{n}\int_{-\pi}^{\pi}\Biggl\vert
\sum_{j=1}^{t-1}\sum_{k=t-j}^{\infty}\frac{\phi_{k}(
\bbvarphi_{0}) }{j}e^{i( j+k) \mu}\Biggr\vert
^{2}\,d\mu\\
&&\qquad\leq\frac{K}{n}\dsum_{t=2}^{n}\dsum_{j=1}^{t-1}%
\dsum_{k=t-j}^{\infty}\dsum_{l=1}^{t-1}\frac{\phi_{k}(%
\bbvarphi_{0})\phi_{j+k-l}(\bbvarphi_{0})}{jl} \\
&&\qquad\leq\frac{K}{n}\dsum_{t=2}^{n}\dsum_{j=1}^{t-1}\dsum%
_{k=t-j}^{\infty}\dsum_{l=1}^{t-1}\frac{k^{-1-\varsigma
}(j+k-l)^{-1-\varsigma}}{jl} \\
&&\qquad\leq\frac{K}{n}\dsum_{t=2}^{n}\dsum_{l=1}^{t-1}\frac{%
(t-l)^{-1-\varsigma}}{l}\dsum_{j=1}^{t-1}\frac{(t-j)^{-\varsigma}}{j}.
\end{eqnarray*}
Now
\begin{eqnarray*}
\sum_{l=1}^{t-1}\frac{( t-l) ^{-1-\varsigma}}{l}%
&=&\sum_{l=1}^{[ t/2] }\frac{( t-l) ^{-1-\varsigma}}{l%
}+\sum_{l=[ t/2] +1}^{t-1}\frac{( t-l)
^{-1-\varsigma}}{l}\\
&\leq& K( t^{-1-\varsigma}\log t+t^{-1}) \leq\frac
{K}{t}.
\end{eqnarray*}
Then $\operatorname{Var}( r_{2}) =O(
n^{-1}\dsum_{t=2}^{n}t^{-1}\dsum_{j=1}^{t-1}j^{-1})
=O( n^{-1}\log^{2}n) $, so
\[
r_{2}=O_{p}(n^{-1/2}\log n).
\]
Next, by Lemma \ref{lemma2}
\[
r_{3}=O_{p}\Biggl( n^{-{1/2}}\sum_{t=2}^{n}t^{-{1/2}-\varsigma
}\log t\Biggr) =O_{p}( n^{-{1/2}}) .
\]
Also, $E( \mathbf{s}_{1}) =0$ and
\begin{eqnarray*}
\operatorname{Var}( \mathbf{s}_{1}) &=&O\Biggl( \Biggl\Vert\frac{1}{n}%
\sum_{t=2}^{n}\sum_{j=t}^{\infty}\sum_{k=t}^{\infty}%
\frac{\partial\phi_{j}( \bbvarphi_{0}) }{\partial
\bbvarphi}\,\frac{\partial\phi_{k}( \bbvarphi%
_{0}) }{\partial\bbvarphi^{\prime}}E(
u_{t-j}u_{t-k}) \Biggr\Vert\Biggr) \\
&=&O\Biggl( \frac{1}{n}\sum_{t=2}^{n}\dint_{-\pi}^{\pi
}\Biggl\Vert\sum_{j=t}^{\infty}\frac{\partial\phi_{j}(
\bbvarphi_{0}) }{\partial\bbvarphi}e^{ij\lambda
}\Biggr\Vert^{2}\,d\lambda\Biggr) \\
&=&O\Biggl( \frac{1}{n}\sum_{t=2}^{n}\sum_{j=t}^{\infty
}\biggl\Vert\frac{\partial\phi_{j}( \bbvarphi_{0}) }{%
\partial\bbvarphi}\biggr\Vert^{2}\Biggr) =O\Biggl( \frac{1}{n}%
\sum_{t=2}^{n}t^{-1-2\varsigma}\Biggr) =O( n^{-1}) ,
\end{eqnarray*}
since $\varsigma>\frac{1}{2}$, \mbox{$\Vert\cdot\Vert$} denoting Euclidean
norm. Finally, by Lemmas \ref{lemma2} and \ref{lemma4}
\[
\mathbf{s}_{2}=O_{p}\Biggl( n^{-{1/2}}\sum_{t=1}^{n}t^{-{1/2%
}-\varsigma}\Biggr) =O_{p}(n^{-{1/2}}),
\]
to conclude the proof of (\ref{x4}).

Next, (\ref{x5}) holds by the Cram\'er--Wold device and, for example,
Theorem~1 of \cite{brown} on showing that
%
%e2.57 #&#
\begin{equation} \label{x7}
E\Biggl( \varepsilon_{t}\sum_{j=1}^{\infty}\mathbf{m}%
_{j}( \bbvarphi_{0}) \varepsilon_{t-j}\Big\vert
\mathcal{F}_{t-1}\Biggr) =0\qquad\mbox{a.s.}
\end{equation}
and
%
%e2.58 #&#
\begin{eqnarray} \label{x56}
&&\frac{1}{n}\sum_{t=2}^{n}E\Biggl( \varepsilon
_{t}^{2}\sum_{j=1}^{\infty}\sum_{k=1}^{\infty}\mathbf{m}%
_{j}( \bbvarphi_{0}) \mathbf{m}_{k}^{\prime}(
\bbvarphi_{0}) \varepsilon_{t-j}\varepsilon_{t-k}\Big\vert
\mathcal{F}_{t-1}\Biggr)
\nonumber\\[-8pt]\\[-8pt]
&&\qquad{}
-\frac{1}{n}\sum_{t=2}^{n}E\Biggl( \varepsilon
_{t}^{2}\sum_{j=1}^{\infty}\sum_{k=1}^{\infty}\mathbf{m}%
_{j}( \bbvarphi_{0}) \mathbf{m}_{k}^{\prime}(
\bbvarphi_{0}) \varepsilon_{t-j}\varepsilon_{t-k}\Biggr)
\rightarrow_{p}0,\nonumber
\end{eqnarray}
because $E( \varepsilon_{t}^{2}\sum_{j=1}^{\infty
}\sum_{k=1}^{\infty}\mathbf{m}_{j}( \bbvarphi%
_{0}) \mathbf{m}_{k}^{\prime}( \bbvarphi_{0})
\varepsilon_{t-j}\varepsilon_{t-k}\vert\mathcal{F}_{t-1}) $
has expectation $\sigma_{0}^{2}\mathbf{A}$, noting that the Lindeberg
condition is satisfied as $\varepsilon_{t}\sum_{j=1}^{\infty}%
\mathbf{m}_{j}( \bbvarphi_{0}) \varepsilon_{t-j}$ is
stationary with finite variance. Now (\ref{x7}) follows as $\varepsilon
_{t-j}$, $j\geq1$, is $\mathcal{F}_{t-1}$-measurable, whereas the
left-hand side
of (\ref{x56}) is
\[
\frac{\sigma_{0}^{2}}{n}\sum_{t=2}^{n}\sum_{j=1}^{\infty
}\sum_{k=1}^{\infty}\mathbf{m}_{j}( \bbvarphi%
_{0}) \mathbf{m}_{k}^{\prime}( \bbvarphi_{0})
\bigl( \varepsilon_{t-j}\varepsilon_{t-k}-E( \varepsilon
_{t-j}\varepsilon_{t-k}) \bigr) \rightarrow_{p}0,
\]
because $\sum_{j=1}^{\infty}\sum_{k=1}^{\infty}\mathbf{m}%
_{j}( \bbvarphi_{0}) \mathbf{m}_{k}^{\prime}(
\bbvarphi_{0}) ( \varepsilon_{t-j}\varepsilon
_{t-k}-E( \varepsilon_{t-j}\varepsilon_{t-k}) ) $ is
stationary ergodic with mean zero. This completes the proof of (\ref
{x5}), and thus (\ref%
{x2}).

\textit{Proof of} (\ref{x3}). Denote by $N_{\epsilon}$ an open
neighborhood of radius $\epsilon<1/2$ about $\bbtau_{0}$, and
%
%e2.59 #&#
%e2.60 #&#
\begin{eqnarray}
\label{x12}
\mathbf{A}_{n}( \bbtau) &=&\frac{1}{n}\sum%
_{t=2}^{n} \sum_{j=0}^{t-1}\sum_{k=1}^{t-1}\biggl(c_{j}%
\,\frac{\partial^{2}c_{k}}{\partial\bbtau\,\partial\bbtau
^{\prime}}+\frac{\partial c_{j}%
}{\partial\bbtau}\,\frac{\partial c_{k}}{\partial\bbtau%
^{\prime}}\biggr) \gamma_{k-j},\\
\label{h12}
\mathbf{A}( \bbtau) &=&\sum_{j=0}^{\infty
}\sum_{k=1}^{\infty}\biggl(c_{j}%
\,\frac{\partial^{2}c_{k}}{\partial\bbtau\,\partial\bbtau
^{\prime}}+\frac{\partial c_{j}%
}{\partial\bbtau}\,\frac{\partial c_{k}}{\partial\bbtau%
^{\prime}}\biggr) \gamma_{k-j}.
\end{eqnarray}
Trivially,%
\[
\frac{1}{2}\,\frac{\partial^{2}R_{n}( \overline{\bbtau}) }{%
\partial\bbtau\,\partial\bbtau^{\prime}}=\frac
{1}{2}\,\frac{%
\partial^{2}R_{n}( \overline{\bbtau}) }{\partial\bbtau\,\partial\bbtau^{\prime}}-\mathbf{A}_{n}( \overline{%
\bbtau}) +\mathbf{A}_{n}( \overline{\bbtau}%
) -\mathbf{A}( \overline{\bbtau}) +\mathbf{A}(
\overline{\bbtau}) -\mathbf{A}( \bbtau_{0})
+\mathbf{A}( \bbtau_{0}) .
\]
Because $c_{j}(\bbtau_{0})=\phi_{j}(\bbtau_{0})$%
, it follows that $\dsum_{j=0}^{\infty}c_{j}(\bbtau%
_{0})u_{t-j}=\varepsilon_{t}$, so the first term in $\mathbf{A}(
\bbtau_{0}) $ is identically zero. Also, as in the proof of (%
\ref{x5}), the second term of $\mathbf{A}( \bbtau_{0}) $
is identically $\sigma_{0}^{2}\mathbf{A}$. Thus, given that by Slutzky's
theorem and continuity of $\mathbf{A}( \bbtau) $ at $%
\bbtau_{0}$, $\mathbf{A}( \overline{\bbtau}) -%
\mathbf{A}( \bbtau_{0}) =o_{p}( 1) $, (\ref{x3}%
) holds on showing
%
%e2.61 #&#
%e2.62 #&#
\begin{eqnarray}
\label{x9}
\sup_{\bbtau\in N_{\epsilon}
}\biggl\Vert\frac{1}{2}\,\frac{\partial^{2}R_{n}( \bbtau)
}{\partial\bbtau\,\partial\bbtau^{\prime}}-\mathbf{A}%
_{n}( \bbtau) \biggr\Vert &=&o_{p}( 1) ,
\\
\label{xxx}
{\sup_{\bbtau\in N_{\epsilon}
}}\Vert\mathbf{A}_{n}( \bbtau) -\mathbf{A}(
\bbtau) \Vert &=&o_{p}( 1)
\end{eqnarray}
for some $\epsilon>0$, as $n\rightarrow\infty$. As $\epsilon
<1/2$, the proof for (\ref{x9}) is almost identical to that for (\ref{4}),
noting the orders in Lemma \ref{lemma4}. To prove (\ref{xxx}), we show that%
%
%e2.63 #&#
\begin{equation}\label{h13}
\sup_{\bbtau\in N_{\epsilon}
}\Biggl\Vert\frac{1}{n}\sum_{t=2}^{n}\sum_{j=0}^{t-1}\sum%
_{k=1}^{t-1}c_{j}\,\frac{\partial^{2}c_{k}}{\partial\bbtau%
\,\partial\bbtau^{\prime}}\gamma_{k-j}-\sum_{j=0}^{\infty
}\sum_{k=1}^{\infty}c_{j}\,
\frac{\partial^{2}c_{k}}{\partial\bbtau\,\partial\bbtau^{\prime}}\gamma_{k-j}
\Biggr\Vert
\end{equation}
is $o_{p}( 1) $, the proof for the corresponding result
concerning the difference between the second terms in (\ref{x12}), (\ref
{h12}%
) being almost identical. By Lemma~\ref{lemma4}, (\ref{h13}) is bounded by%
%
%e2.64 #&#
\begin{eqnarray} \label{h15}\quad
&&\frac{K}{n}\sum_{t=1}^{n}\sum_{j=1}^{t}\sum_{k=t+1}^{%
\infty}j^{\epsilon-1}k^{\epsilon-1}( k-j) ^{-1-\varsigma}\log
^{2}k+\frac{K}{n}\sum_{t=1}^{n}\sum_{j=t}^{\infty
}j^{2\epsilon-2}\log^{2}j \nonumber\\[-8pt]\\[-8pt]
&&\qquad{}+\frac{K}{n}\sum_{t=1}^{n}\sum_{j=t}^{\infty
}\sum_{k=j+1}^{\infty}j^{\epsilon-1}k^{\epsilon-1}(
k-j) ^{-1-\varsigma}\log^{2}k,\nonumber
\end{eqnarray}
noting that (\ref{h14}) implies that $\gamma_{j}=O( j^{-1-\varsigma
}) $. The first term in (\ref{h15}) is bounded by%
%
%e2.65 #&#
\begin{equation} \label{h16}\qquad
\frac{K}{n}\sum_{t=1}^{n}t^{\epsilon}\sum_{k=t+1}^{\infty
}k^{\epsilon+a-1}( k-t) ^{-1-\varsigma}\leq\frac{K}{n}%
\sum_{t=1}^{n}t^{\epsilon}\sum_{k=1}^{\infty}(
k+t) ^{\epsilon+a-1}k^{-1-\varsigma}
\end{equation}
for any $a>0$. Choosing $a$ such that $2\epsilon+a<1$, (\ref{h16}) is
bounded by%
\[
\frac{K}{n}\sum_{t=1}^{n}t^{2\epsilon+a-1}\sum_{k=1}^{\infty
}k^{-1-\varsigma}=O( n^{2\epsilon+a-1}) =o( 1) .
\]
Similarly, the second term in (\ref{h15}) can be easily shown to be $o(
1) $, whereas the third term is bounded by%
%
%e2.66 #&#
\begin{equation} \label{h17}
\frac{K}{n}\sum_{t=1}^{n}\sum_{j=t}^{\infty}j^{2\epsilon
+a-2}\sum_{k=j+1}^{\infty}( k-j) ^{-1-\varsigma}
\end{equation}
for any $a>0$, so choosing again $a$ such that $2\epsilon+a<1$, (\ref{h17})
is $O( n^{2\epsilon+a-1}) =o( 1) $, to conclude the
proof of (\ref{x3}), and thus of the theorem.
\end{pf}

%s3 #&#
\section{Multivariate extension}\label{sec3}

When observations on several related time series are available joint
modeling can achieve efficiency gains. We consider a vector $\mathbf{x}%
_{t}=(x_{1t},\ldots,x_{rt})^{\prime}$ given by%
%
%e3.1 #&#
\begin{equation} \label{zc}
\mathbf{x}_{t}=\bbLamda_{0}^{-1}\{ \mathbf{u}_{t}
\mathbh{1}( t>0)
\},\qquad t=0,\pm1,\ldots,
\end{equation}
where $\mathbf{u}_{t}=( u_{1t},\ldots,u_{rt}) ^{\prime}$,%
%
%e3.2 #&#
\begin{equation} \label{zb}
\mathbf{u}_{t}=\bbTheta( L;\bbvarphi_{0}) \bbvarepsilon_{t},\qquad t=0,\pm1,\ldots,
\end{equation}
in which $\bbvarepsilon_{t}=( \varepsilon_{1t},\ldots,\varepsilon
_{rt}) ^{\prime}$, $\bbvarphi_{0}$ is (as in the univariate
case) a $p\times1$ vector of short-memory parameters, $\bbTheta(s;%
\bbvarphi)=\dsum_{j=0}^{\infty}\bbTheta_{j}(%
\bbvarphi)s^{j}$, $\bbTheta_{0}(\bbvarphi)=I_{r}$
for all~$\bbvarphi$, and $\bbLamda_{0}=\operatorname{diag}( \Delta
^{\delta_{01}},\ldots,\Delta^{\delta_{0r}}) $, where the memory
parameters $\delta_{0i}$ are unknown real numbers. In general, they
can all
be distinct but for the sake of parsimony we allow for the possibility that
they are known to lie in a set of dimension $q<r$. For example, perhaps
as a
consequence of pre-testing, we might believe some or all the $\delta_{0i}$
are equal, and imposing this restriction in the estimation could further
improve efficiency. We introduce known functions $\delta_{i}=\delta
_{i}(%
\bbdelta)$, $i=1,\ldots,r$, of $q\times1$ vector~$\bbdelta$,
such that for some $\bbdelta_{0}$ we have $\delta_{0i}=\delta
_{i}(%
\bbdelta_{0})$, $i=1,\ldots,r$. We denote $\bbtau=(\bbdelta^{\prime},\bbvarphi^{\prime})^{\prime}$ and define
[cf. (\ref{d})]
\[
\bbvarepsilon_{t}(\bbtau)=\bbTheta^{-1}(L;\bbvarphi)\bbLamda( \bbdelta)
\mathbf{x}_{t},\qquad
t\geq1,
\]
where $\bbLamda( \bbdelta) =\operatorname{diag}( \Delta
^{\delta_{1}},\ldots,\Delta^{\delta_{r}}) $. Gaussian likelihood
considerations suggest the multivariate analogue to (\ref{f})
%
%e3.3 #&#
\begin{equation} \label{36}
R_{n}^{\ast}(\bbtau)=\det\{ \bbSigma_{n}(\bbtau)\} ,
\end{equation}
where $\bbSigma_{n}(\bbtau)=n^{-1}\dsum_{t=1}^{n}%
\bbvarepsilon_{t}(\bbtau)\bbvarepsilon_{t}^{\prime}(%
\bbtau)$, assuming that no prior restrictions link~$\bbtau%
_{0}$ with the covariance matrix of $\bbvarepsilon_{t}$.
Unfortunately our consistency proof for the univariate case does not
straightforwardly extend to an estimate minimizing (\ref{36}) if $q>1$.
Also (\ref{36}) is liable to pose a more severe computational challenge
than (\ref{f}) since $p$ is liable to be larger in the multivariate
case and~$q$ may exceed 1; it may be difficult to locate an approximate minimum
of~(\ref{36}) as a preliminary to iteration. We avoid both these problems by
taking a single Newton step from an initial $\sqrt{n}$-consistent
estimate~$\widetilde{\bbtau}$. Defining
\begin{eqnarray*}
\mathbf{H}_{n}( \bbtau) &=&\frac{1}{n}\sum%
_{t=1}^{n}\biggl( \frac{\partial\bbvarepsilon_{t}(
\bbtau) }{\partial\bbtau^{\prime}}\biggr) ^{\prime}%
\bbSigma_{n}^{-1}( \bbtau) \,
\frac{\partial\bbvarepsilon_{t}( \bbtau) }{\partial\bbtau%
^{\prime}}, \\
\mathbf{h}_{n}( \bbtau) &=&\frac{1}{n}\sum%
_{t=1}^{n}\biggl( \frac{\partial\bbvarepsilon_{t}(
\bbtau) }{\partial\bbtau^{\prime}}\biggr) ^{\prime}%
\bbSigma_{n}^{-1}( \bbtau) \bbvarepsilon%
_{t}( \bbtau) ,
\end{eqnarray*}
we consider the estimate
%
%e3.4 #&#
\begin{equation}\label{39}
\widehat{\bbtau}=\widetilde{\bbtau}-\mathbf{H}_{n}^{-1}(
\widetilde{\bbtau})\mathbf{h}_{n}(\widetilde{\bbtau}).\vadjust{\goodbreak}
\end{equation}

We collect together all the requirements for asymptotic normality of $%
\widehat{\bbtau}$ in:

\begin{enumerate}[A4.]
\item[A4.]

\begin{enumerate}[(iiv)]
\item[(i)] For all $\bbvarphi$, $\Theta( e^{i\lambda};\bbvarphi) $ is differentiable in $\lambda$ with derivative in $%
\operatorname{Lip}( \varsigma) $, $\varsigma>1/2;$

\item[(ii)] for all $\bbvarphi$, $\det\{ \bbTheta%
( s;\bbvarphi) \} \neq0, \vert s\vert
\leq1;$

\item[(iii)] the $\bbvarepsilon_{t}$ in (\ref{zb}) are stationary
and ergodic with finite fourth moment, $E( \bbvarepsilon%
_{t}\vert\mathcal{F}_{t-1}) =0$, $E( \bbvarepsilon_{t}\bbvarepsilon_{t}^{\prime}\vert\mathcal{F}%
_{t-1}) =\bbSigma_{0}$ almost surely, where $\bbSigma%
_{0}$ is positive definite, $\mathcal{F}_{t}$ is the $%
\sigma$-field of events generated by $\bbvarepsilon_{s}$,
$s\leq t$, and conditional (on $\mathcal{F}_{t-1}$) third and fourth moments
and cross-moments of elements of $\bbvarepsilon_{t}$ equal
the corresponding unconditional moments;

\item[(iv)] for all $\lambda$, $\Theta( e^{i\lambda};\bbvarphi) $ is
twice continuously differentiable in $\bbvarphi$ on
a~closed neighborhood $\mathcal{N}_{\epsilon}( \bbvarphi%
_{0}) $ of radius $0<\epsilon<1/2$ about $\bbvarphi_{0};$

\item[(v)] the matrix $B$ having $( i,j) $th element%
\[
\sum_{k=1}^{\infty}\operatorname{tr}\bigl\{ \bigl( \mathbf{d}_{k}^{( i)
}( \bbvarphi_{0}) \bigr) ^{\prime}\bbSigma%
_{0}^{-1}\mathbf{d}_{k}^{( j) }( \bbvarphi%
_{0}) \bbSigma_{0}\bigr\}
\]
is nonsingular, where%
\begin{eqnarray*}
\mathbf{d}_{k}^{( i) }( \bbvarphi_{0}) &=&-%
\frac{\partial\delta_{i}( \bbdelta_{0}) }{\partial
\delta_{i}}\sum_{l=1}^{k}\frac{1}{l}\sum_{m=0}^{k-l}\bbPhi_{m}^{( i) }( \bbvarphi_{0}) \bbTheta_{k-l-m}( \bbvarphi_{0}) ,
\qquad 1\leq i\leq r, \\
&=&\sum_{l=1}^{k}\frac{\partial\bbPhi_{l}( \bbvarphi_{0}) }{\partial\varphi_{i}}\bbTheta_{k-l}(
\bbvarphi_{0}) ,\qquad r+1\leq i\leq r+p,
\end{eqnarray*}
the $\Phi_{j}( \bbvarphi) $ being coefficients in the
expansion $\bbTheta^{-1}( s;\bbvarphi)\,{=}\,\bbPhi( s,\bbvarphi)\,{=}\allowbreak\dsum_{j=0}^{\infty}%
\bbPhi_{j}( \bbvarphi) s^{j}$, where $\Phi
_{m}^{( i) }( \bbvarphi_{0}) $ is an $r\times
r$ matrix whose $i$th column is the $i$th column of $\Phi_{i}(
\bbvarphi_{0}) $ and whose other elements are all zero;

\item[(vi)] $\delta_{i}( \bbdelta) $ is twice
continuously differentiable in $\bbdelta$, for $i=1,\ldots,r;$

\item[(vii)] $\widetilde{\bbtau}$ is a $\sqrt{n}$-consistent
estimate of $\bbtau_{0}$.
\end{enumerate}
\end{enumerate}

The components of A4 are mostly natural extensions of ones in A1, A2
and~A3,
are equally checkable, and require no additional discussion. The important
exception is (vii). When $\bbTheta(s;\bbvarphi)$ is a
diagonal matrix [as in the simplest case $\bbTheta(s;\bbvarphi)\equiv\mathbf{I}_{r}$, when $x_{it}$ is a FARIMA$(0,\delta_{0i},0)$
for $%
i=1,\ldots,r$] then $\widetilde{\bbtau}$ can be obtained by first
carrying out $r$ univariate fits following the approach of
Section~\ref{sec2}, and then if necessary reducing the dimensionality
in a common-sense way: for example, if some of the $\delta_{0i}$ are a
priori equal then the common memory parameter might be estimated by the
arithmetic mean of estimates from the relevant univariate fits. Notice
that in the diagonal-$\bbTheta$ case with no cross-equation parameter
restrictions the efficiency improvement afforded by $\widehat{\bbtau}$
is due solely to cross-correlation in $\bbvarepsilon_{t}$, that is,
nondiagonality of $\bbSigma_{0}$.

When $\bbTheta(s;\bbvarphi)$ is not diagonal, it is less clear how to
use the $\sqrt{n}$-consistent outcome of Theorem \ref{theo2.2} to form
$%
\widetilde{\bbtau}$. We can infer that $\mathbf{u}_{t}$ has spectral
density matrix $(2\pi)^{-1}\bbTheta(e^{i\lambda};\bbvarphi%
_{0})\bbSigma_{0}\bbTheta(e^{-i\lambda};\bbvarphi%
_{0})^{\prime}$. From the $i$th diagonal element of this (the power
spectrum of $u_{it}$), we can deduce a form for the Wold representation
of $%
u_{it}$, corresponding to (\ref{b}). However, starting from innovations~$\bbvarepsilon_{t}$
in (\ref{zb}) satisfying (iii) of A4, it
does not
follow in general that the innovations in the Wold representation of $u_{it}$
will satisfy a condition analogous to (\ref{28}) of A2, indeed it does not
help if we simply strengthen A4 such that the $\bbvarepsilon_{t}$
are independent and identically distributed. However, (\ref{28}) certainly
holds if $\bbvarepsilon_{t}$ is Gaussian, which motivates our
estimation approach from an efficiency perspective. Notice that if
$\mathbf{u}_{t}$ is a vector ARMA process with nondiagonal $\bbTheta$, in
general all $r$ univariate AR operators are identical, and of possibly high
degree; the formation of $\widetilde{\bbtau}$ is liable to
be affected by a~lack of parsimony, or some ambiguity.

An alternative approach could involve first estimating the $\delta
_{0i}$ by
some semiparametric approach, using these estimates to form differenced
$%
\mathbf{x}_{t}$ and then estimating $\bbvarphi_{0}$ from these
proxies for $\mathbf{u}_{t}$. This initial estimate will be
less-than-$\sqrt{%
n}$-consistent, but its rate can be calculated given a rate for the
bandwidth used in the semiparametric estimation. One can then calculate the
(finite) number of iterations of form (\ref{39}) needed to produce an
estimate satisfying (\ref{213}), following Theorem 5 and the discussion on
page 539 of \cite{robinson1}.

%th3.1 #&#
\begin{theorem}
Let (\ref{zc}), (\ref{zb}) and \textup{A4} hold. Then as $n\rightarrow
\infty$%
%
%e3.5 #&#
\begin{equation} \label{zh}
n^{{1/2}}( \widehat{\bbtau}-\bbtau_{0})
\rightarrow_{d}N( \mathbf{0},\mathbf{B}^{-1}) .
\end{equation}
\end{theorem}
\begin{pf}
Because $\widehat{\bbtau}$ is explicitly defined in (\ref{39}), we
start, standardly, by approximating $h_{n}( \widetilde{\bbtau}%
) $ by the mean value theorem. Then in view of A4(vii),
(\ref{zh})~follows on showing
%
%e3.6 #&#
\begin{eqnarray} \label{x15}
\sqrt{n}\mathbf{h}_{n}(\bbtau_{0})&\rightarrow_{d}&N( \mathbf{0},
\mathbf{B}) ,
\\
%
%
%e3.7 #&#
\label{x16}
\mathbf{H}_{n}(\bbtau_{0})&\rightarrow_{p}&\mathbf{B},
\\
%
%
%e3.8 #&#
\label{x17}
\mathbf{H}_{n}( \overline{\bbtau}) -\mathbf{H}_{n}(
\bbtau_{0}) &\rightarrow_{p}&0
\end{eqnarray}
for $\Vert\overline{\bbtau}-\bbtau_{0}\Vert\leq
\Vert\widetilde{\bbtau}-\bbtau_{0}\Vert$. We
only show (\ref{x15}), as (\ref{x16}), (\ref{x17}) follow from similar
arguments to those given in the proof of (\ref{x3}). Noting that
$\partial
\bbvarepsilon_{1}(\bbtau_{0})/\allowbreak\partial\bbtau%
^{\prime}=0$, whereas for $t\geq2$, $\partial\bbvarepsilon_{t}(%
\bbtau_{0})/\partial\bbtau^{\prime}$ equals
\begin{eqnarray*}
&&\sum_{j=1}^{t-1}\Biggl( -\bbPhi_{j}^{( 1) }(
\bbvarphi_{0}) \sum_{k=1}^{t-j-1}\frac{\mathbf{u}%
_{t-j-k}}{k},\ldots,-\bbPhi_{j}^{( r) }( \bbvarphi_{0}) \sum_{k=1}^{t-j-1}\frac{\mathbf{u}_{t-j-k}}{k}%
, \\
&&\hspace*{120.5pt} \frac{\partial\bbPhi_{j}( \bbvarphi%
_{0}) }{\partial\varphi_{1}}\mathbf{u}_{t-j},\ldots,\frac{\partial
\bbPhi_{j}( \bbvarphi_{0}) }{\partial\varphi_{p}%
}\mathbf{u}_{t-j}\Biggr)
\end{eqnarray*}
by similar arguments to those in the proof of Theorem \ref{theo2.2}, it can be shown
that the left-hand side of (\ref{x15}) equals%
\[
\frac{1}{\sqrt{n}}\sum_{t=2}^{n}
\Biggl(
\sum_{j=1}^{\infty}\mathbf{d}_{j}^{( 1) }( \bbvarphi_{0}) \bbvarepsilon_{t-j} \cdots
\sum_{j=1}^{\infty}\mathbf{d}_{j}^{( r+p) }( \bbvarphi_{0}) \bbvarepsilon_{t-j}%
\Biggr) ^{\prime}\bbSigma_{0}^{-1}\bbvarepsilon%
_{t}+o_{p}( 1) .
\]
Then by the Cram\'er--Wold device, (\ref{x15}) holds if for any $(
r+p) $-dimensional vector~$\bbvartheta$ (with $i$th component
$\vartheta_{i}$)
%
%e3.9 #&#
\begin{equation}\label{x20}
\frac{1}{\sqrt{n}}\sum_{t=2}^{n}\sum_{j=1}^{\infty}\bbvarepsilon_{t-j}^{\prime}\mathbf{M}_{j}^{\prime}( \bbvarphi
_{0}) \bbSigma_{0}^{-1}\bbvarepsilon_{t}\rightarrow
_{d}N( 0,\bbvartheta^{\prime}\mathbf{B}\bbvartheta) ,
\end{equation}
where $\mathbf{M}_{j}( \bbvarphi_{0})
=\dsum_{k=1}^{r+p}\vartheta_{k}\mathbf{d}_{j}^{( k)
}( \bbvarphi_{0})$. As in the proof of (\ref{x5}), (\ref{x20})
holds by Theorem 1 of
\cite{brown}, for example, noting that
\begin{eqnarray*}
&&E\Biggl( \sum_{j=1}^{\infty}\bbvarepsilon_{t-j}^{\prime}%
\mathbf{M}_{j}^{\prime}( \bbvarphi_{0}) \bbSigma%
_{0}^{-1}\bbvarepsilon_{t}\Biggr) ^{2} \\
&&\qquad=E\Biggl( \sum_{j=1}^{\infty}\sum_{k=1}^{\infty}\bbvarepsilon_{t-j}^{\prime}\mathbf{M}_{j}^{\prime}( \bbvarphi
_{0}) \bbSigma_{0}^{-1}E( \bbvarepsilon_{t}%
\bbvarepsilon_{t}^{\prime}\vert\mathcal{F}_{t-1})
\bbSigma_{0}^{-1}\mathbf{M}_{k}( \bbvarphi_{0})
\bbvarepsilon_{t-k}\Biggr) \\
&&\qquad=E\Biggl( \sum_{j=1}^{\infty}\sum_{k=1}^{\infty}\operatorname{tr}\{
\bbvarepsilon_{t-j}^{\prime}\mathbf{M}_{j}^{\prime}( \bbvarphi_{0}) \bbSigma_{0}^{-1}\mathbf{M}_{k}( \bbvarphi_{0}) \bbvarepsilon_{t-k}\} \Biggr) \\
&&\qquad=\sum_{j=1}^{\infty}\operatorname{tr}\{ \mathbf{M}_{j}^{\prime}(
\bbvarphi_{0}) \bbSigma_{0}^{-1}\mathbf{M}_{j}(
\bbvarphi_{0}) \bbSigma_{0}\} =\bbvartheta^{\prime}\mathbf{B}\bbvartheta
\end{eqnarray*}
to conclude the proof.
\end{pf}

%s4 #&#
\section{Further comments and extensions}\label{sec4}

(1) Our univariate and multivariate structures cover a wide range of
parametric models for stationary and nonstationary time series, with memory
parameters allowed to lie in a set that can be arbitrarily large. Unit root
series are a special case, but unlike in the bulk of the large unit root
literature, we do not have to assume knowledge that memory parameters
are 1.
Indeed, in Monte Carlo \cite{hualde} our method out-performs one which
correctly assumes the unit interval in which~$\delta_{0}$ lies, while in
empirical examples our findings conflict with previous, unit root, ones.

(2) As the nondiagonal structure of $\mathbf{A}$ and $\mathbf{B}$
suggests, there is efficiency loss in estimating $\bbvarphi_{0}$ if
memory parameters are unknown, but on the other hand if these are
misspecified, $\bbvarphi_{0}$ will in general be inconsistently
estimated. Our limit distribution theory can be used to test hypotheses on
the memory and other parameters, after straightforwardly forming consistent
estimates of $\mathbf{A}$ or $\mathbf{B}$.

(3) Our multivariate system (\ref{zc}), (\ref{zb}) does not cover
fractionally cointegrated systems because $\bbSigma_{0}$
is required to be positive definite. On the other hand, our theory for
univariate estimation should cover estimation of individual memory
parameters, so long as Assumption A2, in particular, can be reconciled with
the full system specification. Moreover, again on an individual basis, it
should be possible to derive analogous properties of estimates of memory
parameters of cointegrating errors based on residuals that use simple
estimates of cointegrating vectors, such as least squares.

(4) In a more standard regression setting, for example, with
deterministic regressors such as polynomial functions of time, it
should be
possible to extend our theory for univariate and multivariate models to
residual-based estimates of memory parameters of errors.

(5) Adaptive estimates, which have greater efficiency at distributions
of unknown, non-Gaussian form, can be obtained by taking one Newton step
from our estimates (as in \cite{robinson2}).

(6) Our methods of proof should be extendable to cover seasonally and
cyclically fractionally differenced processes.

(7) Nonstationary fractional series can be defined in many ways. Our
definition [(\ref{a}) and (\ref{zc})] is a leading one in the literature,
and has been termed ``Type II.'' Another
popular one (``Type I'') was used by
\cite{velasco} for an alternate type of estimate. That estimate assumes
invertibility and is generally less efficient than $\widehat{\bbtau}$
due to the tapering required to handle nonstationarity. It seems likely that
the asymptotic theory derived in this paper for $\widehat{\bbtau}$
can also be established in a ``Type
I'' setting.

%s5 #&#
\section{Technical lemmas}\label{sec5}

The proofs of the following lemmas appear in \cite{hualde}.

%le1 #&#
\begin{lemma}\label{lemma1}
Under \textup{A1}
%
%e5.1 #&#
\begin{equation}\label{aad}
\varepsilon_{t}( \bbtau)
=\sum_{j=0}^{t-1}c_{j}( \bbtau) u_{t-j}
\end{equation}
with $c_{0}( \bbtau) =1$ where for any %
$\delta\in\mathcal{I}$, as $j\rightarrow\infty$,
%
%e5.2 #&#
\begin{eqnarray}\label{x21}
\sup_{\bbvarphi\in\Psi}\vert c_{j}( \bbtau%
) \vert&=&O\bigl( j^{\max( \delta_{0}-\delta-1,-1-\varsigma
) }\bigr) ,\nonumber\\[-8pt]\\[-8pt]
{\sup_{\bbvarphi\in\Psi}}\vert c_{j+1}( \bbtau%
) -c_{j}( \bbtau) \vert&=&O\bigl( j^{\max
( \delta_{0}-\delta-2,-1-\varsigma) }\bigr).\nonumber
\end{eqnarray}
\end{lemma}
%
%le2 #&#
\begin{lemma}\label{lemma2}
Under \textup{A1, A2}
\[
\varepsilon_{t}( \bbtau^{\ast})
=\sum_{j=0}^{t-1}a_{j}\varepsilon_{t-j}+v_{t}( \delta) ,
\]
where $\bbtau^{\ast}=( \delta,\bbvarphi%
_{0}) $ and for any $\kappa\geq1/2$%
\[
{\sup_{\delta_{0}-\kappa\leq\delta<\delta_{0}-{1/2}+\eta
}}\vert v_{t}( \delta) \vert=O_{p}( t^{\kappa
-1})
\]
and $v_{t}( \delta_{0}) =O_{p}( t^{-1/2-\varsigma}) .
$
\end{lemma}
%
%le3 #&#
\begin{lemma}\label{lemma3} Under
\textup{A1, A2}
%
%e5.3 #&#
\begin{equation} \label{aah}
\sum_{j=1}^{n}I_{\varepsilon( \bbtau) }(
\lambda_{j}) =\sum_{j=1}^{n}\biggl\vert\frac{\theta(
e^{i\lambda_{j}};\bbvarphi_{0}) }{\theta( e^{i\lambda
_{j}};\bbvarphi) }\biggr\vert^{2}I_{\varepsilon(
\bbtau^{\ast}) }( \lambda_{j}) +V_{n}(
\bbtau) ,
\end{equation}
where for any real number $\kappa\geq1/2$
%
%e5.4 #&#
\begin{equation}
{\mathop{\sup_{\delta_{0}-\kappa\leq\delta<\delta
_{0}-{1/2}+\eta}}_{\bbvarphi\in\Psi}}\vert V_{n}(
\bbtau) \vert=O_{p}\bigl( \log^{2}n
\mathbh{1}( \kappa
=1/2) +n^{2\kappa-1}\mathbh{1}( \kappa>1/2) \bigr) .\hspace*{-32pt}
\end{equation}
\end{lemma}
%
%le4 #&#
\begin{lemma}\label{lemma4}
Under \textup{A3}, given an open neighborhood $%
N_{\epsilon}$ of radius $%
\epsilon<1/2$ about $\bbtau_{0}$, as $%
j\rightarrow\infty$,
\begin{eqnarray*}
\sup_{\bbtau\in N_{\epsilon}
}\vert c_{j}( \bbtau) \vert &=&O(
j^{\epsilon-1}), \\
\sup_{\bbtau\in N_{\epsilon} }\biggl\vert\frac
{\partial c_{j}( \bbtau) }{\partial\delta}\biggr\vert&=&O( j^{\epsilon-1}\log
j), \\
\sup_{\bbtau\in N_{\epsilon}
}\vert c_{j+1}( \bbtau) -c_{j}( \bbtau%
) \vert &=&O\bigl( j^{\max( \epsilon-2,-1-\varsigma)
}\bigr) , \\
\sup_{\bbtau\in N_{\epsilon}
}\biggl\vert\frac{\partial}{\partial\delta}\bigl( c_{j+1}( \bbtau) -c_{j}( \bbtau) \bigr)
\biggr\vert
&=&O( j^{-1-\varsigma}+j^{\epsilon-2}\log j) , \\
\sup_{\bbtau\in N_{\epsilon}
}\biggl\vert\frac{\partial^{2}c_{j}( \bbtau) }{\partial
\delta^{2}}\biggr\vert &=&O( j^{\epsilon-1}\log^{2}j),\\
\sup_{\bbtau\in N_{\epsilon}
}\biggl\Vert\frac{\partial c_{j}( \bbtau) }{\partial
\bbvarphi}\biggr\Vert&=&O( j^{\epsilon-1}), \\
\sup_{\bbtau\in N_{\epsilon}
}\biggl\vert\frac{\partial^{2}}{\partial\delta^{2}}\bigl( c_{j+1}(
\bbtau) -c_{j}( \bbtau) \bigr) \biggr\vert
&=&O( j^{-1-\varsigma}+j^{\epsilon-2}\log^{2}j) ,\\
\sup_{\bbtau\in N_{\epsilon}
}\biggl\Vert\frac{\partial}{\partial\bbvarphi}\bigl( c_{j+1}(
\bbtau) -c_{j}( \bbtau) \bigr) \biggr\Vert
&=&O\bigl( j^{\max( \epsilon-2,-1-\varsigma) }\bigr) , \\
\sup_{\bbtau\in N_{\epsilon}
}\biggl\Vert\frac{\partial^{2}c_{j}( \bbtau) }{\partial
\bbvarphi\,\partial\bbvarphi^{\prime}}\biggr\Vert
&=&O( j^{\epsilon-1}), \\
\sup_{\bbtau\in N_{\epsilon} }\biggl\Vert\frac{\partial
^{2}c_{j}( \bbtau) }{\partial\bbvarphi\,\partial
\delta}\biggr\Vert&=&O( j^{\epsilon-1}\log j), \\
\sup_{\bbtau\in N_{\epsilon}
}\biggl\Vert\frac{\partial^{2}}{\partial\bbvarphi\,\partial\bbvarphi^{\prime}}\bigl( c_{j+1}( \bbtau) -c_{j}(
\bbtau) \bigr)\biggr\Vert &=&O\bigl( j^{\max(
\epsilon-2,-1-\varsigma) }\bigr) , \\
\sup_{\bbtau\in N_{\epsilon}
}\biggl\Vert\frac{\partial^{2}}{\partial\bbvarphi\,\partial\delta}%
\bigl( c_{j+1}( \bbtau) -c_{j}( \bbtau%
) \bigr) \biggr\Vert &=&O( j^{-1-\varsigma}+j^{\epsilon-2}\log
j) .
\end{eqnarray*}
\end{lemma}

\section*{Acknowledgments}

We thank the Associate Editor and two referees for constructive
comments that have improved the presentation. We also thank S\o ren
Johansen and Morten O. Nielsen for helpful comments. Some of the second
author's work was carried out while visiting Universidad Carlos III,
Madrid, holding a C\'{a}tedra de Excelencia.

\begin{supplement}[id=suppA]
\stitle{Supplement to ``Gaussian pseudo-maximum likelihood
estimation of fractional time series models''}
\slink[doi]{10.1214/11-AOS931SUPP} %[doi,text={...}] - jei reikia
%suskaldyti doi
\sdatatype{.pdf}
\sfilename{aos931\_supp.pdf}
\sdescription{The supplementary material contains a Monte Carlo
experiment of finite sample performance of the proposed procedure, an
empirical application to U.S. income and consumption data, and the
proofs of the lemmas given in Section~\ref{sec5} of the present paper.}
\end{supplement}

% imsref loaded by lrinkeviciute, 2011-12-28 12:41:17
%

\printaddresses


\begin{thebibliography}{28}
% BibTex style file: ims.bst, 2011-05-30
% Default style options (sort=0,type=number).
% Used options (sort=1,type=number).

%b1 #&#
\bibitem{abramowitz}
%
\begin{bbook}[auto:STB|2011/12/23|09:17:17]
\bauthor{\bsnm{Abramowitz},~\bfnm{M.}\binits{M.}} \AND
\bauthor{\bsnm{Stegun},~\bfnm{I.}\binits{I.}}
(\byear{1970}).
\btitle{Handbook of Mathematical Functions}.
\bpublisher{Dover}, \baddress{New York}.
\bptok{imsref}%
\end{bbook}
%
\endbibitem

%b2 #&#
\bibitem{adensted}
%
\begin{barticle}[mr]
\bauthor{\bsnm{Adenstedt},~\bfnm{Rolf~K.}\binits{R.~K.}}
(\byear{1974}).
\btitle{On large-sample estimation for the mean of a stationary random
sequence}.
\bjournal{Ann. Statist.}
\bvolume{2}
\bpages{1095--1107}.
\bid{issn={0090-5364}, mr={0368354}}
\bptok{imsref}%
\end{barticle}
%
\endbibitem

%b3 #&#
\bibitem{beran}
%
\begin{barticle}[mr]
\bauthor{\bsnm{Beran},~\bfnm{Jan}\binits{J.}}
(\byear{1995}).
\btitle{Maximum likelihood estimation of the differencing parameter for
invertible short and long memory autoregressive integrated moving average
models}.
\bjournal{J. Roy. Statist. Soc. Ser. B}
\bvolume{57}
\bpages{659--672}.
\bid{issn={0035-9246}, mr={1354073}}
\bptok{imsref}%
\end{barticle}
%
\endbibitem

%b4 #&#
\bibitem{billingsley}
%
\begin{bbook}[mr]
\bauthor{\bsnm{Billingsley},~\bfnm{Patrick}\binits{P.}}
(\byear{1968}).
\btitle{Convergence of Probability Measures}.
\bpublisher{Wiley}, \baddress{New York}.
\bid{mr={0233396}}
\bptok{imsref}%
\end{bbook}
%
\endbibitem

%b5 #&#
\bibitem{bloomfield}
%
\begin{barticle}[mr]
\bauthor{\bsnm{Bloomfield},~\bfnm{P.}\binits{P.}}
(\byear{1973}).
\btitle{An exponential model for the spectrum of a scalar time series}.
\bjournal{Biometrika}
\bvolume{60}
\bpages{217--226}.
\bid{issn={0006-3444}, mr={0323048}}
\bptok{imsref}%
\end{barticle}
%
\endbibitem

%b6 #&#
\bibitem{box}
%
\begin{bbook}[mr]
\bauthor{\bsnm{Box},~\bfnm{George E.~P.}\binits{G.~E.~P.}} \AND
\bauthor{\bsnm{Jenkins},~\bfnm{Gwilym~M.}\binits{G.~M.}}
(\byear{1971}).
\btitle{Time Series Analysis: Forecasting and Control}. %,
\bpublisher{Holden-Day}, \baddress{San Francisco, CA}.
\bptnote{check year}%
\bptok{imsref}%
\end{bbook}
%
\endbibitem

%b7 #&#
\bibitem{brown}
%
\begin{barticle}[mr]
\bauthor{\bsnm{Brown},~\bfnm{B.~M.}\binits{B.~M.}}
(\byear{1971}).
\btitle{Martingale central limit theorems}.
\bjournal{Ann. Math. Statist.}
\bvolume{42}
\bpages{59--66}.
\bid{issn={0003-4851}, mr={0290428}}
\bptok{imsref}%
\end{barticle}
%
\endbibitem

%b8 #&#
\bibitem{dahlhaus}
%
\begin{barticle}[mr]
\bauthor{\bsnm{Dahlhaus},~\bfnm{Rainer}\binits{R.}}
(\byear{1989}).
\btitle{Efficient parameter estimation for self-similar processes}.
\bjournal{Ann. Statist.}
\bvolume{17}
\bpages{1749--1766}.
\bid{doi={10.1214/aos/1176347393}, issn={0090-5364}, mr={1026311}}
\bptok{imsref}%
\end{barticle}
%
\endbibitem

%b9 #&#
\bibitem{fox}
%
\begin{barticle}[mr]
\bauthor{\bsnm{Fox},~\bfnm{Robert}\binits{R.}} \AND
\bauthor{\bsnm{Taqqu},~\bfnm{Murad~S.}\binits{M.~S.}}
(\byear{1986}).
\btitle{Large-sample properties of parameter estimates for strongly dependent
stationary {G}aussian time series}.
\bjournal{Ann. Statist.}
\bvolume{14}
\bpages{517--532}.
\bid{doi={10.1214/aos/1176349936}, issn={0090-5364}, mr={0840512}}
\bptok{imsref}%
\end{barticle}
%
\endbibitem

%b10 #&#
\bibitem{giraitis}
%
\begin{barticle}[mr]
\bauthor{\bsnm{Giraitis},~\bfnm{L.}\binits{L.}} \AND
\bauthor{\bsnm{Surgailis},~\bfnm{D.}\binits{D.}}
(\byear{1990}).
\btitle{A central limit theorem for quadratic forms in strongly dependent
linear variables and its application to asymptotical normality of {W}hittle's
estimate}.
\bjournal{Probab. Theory Related Fields}
\bvolume{86}
\bpages{87--104}.
\bid{doi={10.1007/BF01207515}, issn={0178-8051}, mr={1061950}}
\bptok{imsref}%
\end{barticle}
%
\endbibitem

%b11 #&#
\bibitem{hannan}
%
\begin{barticle}[mr]
\bauthor{\bsnm{Hannan},~\bfnm{E.~J.}\binits{E.~J.}}
(\byear{1973}).
\btitle{The asymptotic theory of linear time-series models}.
\bjournal{J. Appl. Probab.}
\bvolume{10}
\bpages{130--145}.
\bid{issn={0021-9002}, mr={0365960}}
\bptnote{check related}%
\bptok{imsref}%
\end{barticle}
%
\endbibitem

%b12 #&#
\bibitem{hosoya1}
%
\begin{barticle}[mr]
\bauthor{\bsnm{Hosoya},~\bfnm{Yuzo}\binits{Y.}}
(\byear{1996}).
\btitle{The quasi-likelihood approach to statistical inference on multiple
time-series with long-range dependence}.
\bjournal{J. Econometrics}
\bvolume{73}
\bpages{217--236}.
\bid{doi={10.1016/0304-4076(95)01738-0}, issn={0304-4076}, mr={1410005}}
\bptok{imsref}%
\end{barticle}
%
\endbibitem

%b13 #&#
\bibitem{hosoya2}
%
\begin{barticle}[mr]
\bauthor{\bsnm{Hosoya},~\bfnm{Yuzo}\binits{Y.}}
(\byear{2005}).
\btitle{Fractional invariance principle}.
\bjournal{J. Time Ser. Anal.}
\bvolume{26}
\bpages{463--486}.
\bid{doi={10.1111/j.1467-9892.2004.00411.x}, issn={0143-9782}, mr={2163292}}
\bptok{imsref}%
\end{barticle}
%
\endbibitem

%b14 #&#
\bibitem{hualde}
%
\begin{bmisc}[auto:STB|2011/12/23|09:17:17]
\bauthor{\bsnm{Hualde},~\bfnm{J.}\binits{J.}} \AND
\bauthor{\bsnm{Robinson},~\bfnm{P.~M.}\binits{P.~M.}}
(\byear{2011}).
\bhowpublished{Supplement to `` Gaussian pseudo-maximum likelihood
estimation of fractional time series models.''
\href{http://dx.doi.org/10.1214/11-AOS931SUPP}{DOI:10.1214/}
\href{http://dx.doi.org/10.1214/11-AOS931SUPP}{11-AOS931SUPP}.}
\bptok{imsref}%
\end{bmisc}
%
\endbibitem

%b15 #&#
\bibitem{li}
%
\begin{barticle}[mr]
\bauthor{\bsnm{Li},~\bfnm{W.~K.}\binits{W.~K.}} \AND
\bauthor{\bsnm{McLeod},~\bfnm{A.~I.}\binits{A.~I.}}
(\byear{1986}).
\btitle{Fractional time series modelling}.
\bjournal{Biometrika}
\bvolume{73}
\bpages{217--221}.
\bid{doi={10.1093/biomet/73.1.217}, issn={0006-3444}, mr={0836451}}
\bptok{imsref}%
\end{barticle}
%
\endbibitem

%b16 #&#
\bibitem{nordman}
%
\begin{barticle}[mr]
\bauthor{\bsnm{Nordman},~\bfnm{Daniel~J.}\binits{D.~J.}} \AND
\bauthor{\bsnm{Lahiri},~\bfnm{Soumendra~N.}\binits{S.~N.}}
(\byear{2006}).
\btitle{A frequency domain empirical likelihood for short- and long-range
dependence}.
\bjournal{Ann. Statist.}
\bvolume{34}
\bpages{3019--3050}.
\bid{doi={10.1214/009053606000000902}, issn={0090-5364}, mr={2329476}}
\bptok{imsref}%
\end{barticle}
%
\endbibitem

%b17 #&#
\bibitem{robinson1}
%
\begin{barticle}[mr]
\bauthor{\bsnm{Robinson},~\bfnm{P.~M.}\binits{P.~M.}}
(\byear{1988}).
\btitle{The stochastic difference between econometric statistics}.
\bjournal{Econometrica}
\bvolume{56}
\bpages{531--548}.
\bid{doi={10.2307/1911699}, issn={0012-9682}, mr={0946120}}
\bptok{imsref}%
\end{barticle}
%
\endbibitem

%b18 #&#
\bibitem{robinson1a}
%
\begin{barticle}[mr]
\bauthor{\bsnm{Robinson},~\bfnm{P.~M.}\binits{P.~M.}}
(\byear{1994}).
\btitle{Efficient tests of nonstationary hypotheses}.
\bjournal{J. Amer. Statist. Assoc.}
\bvolume{89}
\bpages{1420--1437}.
\bid{issn={0162-1459}, mr={1310232}}
\bptok{imsref}%
\end{barticle}
%
\endbibitem

%b19 #&#
\bibitem{robinson1b}
%
\begin{barticle}[mr]
\bauthor{\bsnm{Robinson},~\bfnm{P.~M.}\binits{P.~M.}}
(\byear{1995}).
\btitle{Gaussian semiparametric estimation of long range dependence}.
\bjournal{Ann. Statist.}
\bvolume{23}
\bpages{1630--1661}.
\bid{doi={10.1214/aos/1176324317}, issn={0090-5364}, mr={1370301}}
\bptok{imsref}%
\end{barticle}
%
\endbibitem

%b20 #&#
\bibitem{robinson2}
%
\begin{barticle}[mr]
\bauthor{\bsnm{Robinson},~\bfnm{P.~M.}\binits{P.~M.}}
(\byear{2005}).
\btitle{Efficiency improvements in inference on stationary and nonstationary
fractional time series}.
\bjournal{Ann. Statist.}
\bvolume{33}
\bpages{1800--1842}.
\bid{doi={10.1214/009053605000000354}, issn={0090-5364}, mr={2166563}}
\bptok{imsref}%
\end{barticle}
%
\endbibitem

%b21 #&#
\bibitem{robinson3}
%
\begin{bincollection}[mr]
\bauthor{\bsnm{Robinson},~\bfnm{P.~M.}\binits{P.~M.}}
(\byear{2006}).
\btitle{Conditional-sum-of-squares estimation of models for stationary time
series with long memory}.
In \bbooktitle{Time Series and Related Topics: In Memory of
Ching-Zong Wei}
(\beditor{H.-C. Ho}, \beditor{C.-K. Ing} and \beditor{T. L. Lai}, eds.).
\bseries{Institute of Mathematical Statistics Lecture Notes---Monograph Series}
\bvolume{52}
\bpages{130--137}.
\bpublisher{IMS}, \baddress{Beachwood, OH}.
\bid{doi={10.1214/074921706000000996}, mr={2427843}}
\bptok{imsref}%
\end{bincollection}
%
\endbibitem

%b22 #&#
\bibitem{robhualde}
%
\begin{barticle}[mr]
\bauthor{\bsnm{Robinson},~\bfnm{P.~M.}\binits{P.~M.}} \AND
\bauthor{\bsnm{Hualde},~\bfnm{J.}\binits{J.}}
(\byear{2003}).
\btitle{Cointegration in fractional systems with unknown integration orders}.
\bjournal{Econometrica}
\bvolume{71}
\bpages{1727--1766}.
\bid{doi={10.1111/1468-0262.00468}, issn={0012-9682}, mr={2015418}}
\bptok{imsref}%
\end{barticle}
%
\endbibitem

%b23 #&#
\bibitem{shimotsu}
%
\begin{barticle}[mr]
\bauthor{\bsnm{Shimotsu},~\bfnm{Katsumi}\binits{K.}} \AND
\bauthor{\bsnm{Phillips},~\bfnm{Peter C.~B.}\binits{P.~C.~B.}}
(\byear{2005}).
\btitle{Exact local {W}hittle estimation of fractional integration}.
\bjournal{Ann. Statist.}
\bvolume{33}
\bpages{1890--1933}.
\bid{doi={10.1214/009053605000000309}, issn={0090-5364}, mr={2166565}}
\bptok{imsref}%
\end{barticle}
%
\endbibitem

%b24 #&#
\bibitem{tanaka}
%
\begin{barticle}[mr]
\bauthor{\bsnm{Tanaka},~\bfnm{Katsuto}\binits{K.}}
(\byear{1999}).
\btitle{The nonstationary fractional unit root}.
\bjournal{Econometric Theory}
\bvolume{15}
\bpages{549--582}.
\bid{doi={10.1017/S0266466699154045}, issn={0266-4666}, mr={1717967}}
\bptok{imsref}%
\end{barticle}
%
\endbibitem

%b25 #&#
\bibitem{velasco}
%
\begin{barticle}[mr]
\bauthor{\bsnm{Velasco},~\bfnm{Carlos}\binits{C.}} \AND
\bauthor{\bsnm{Robinson},~\bfnm{Peter~M.}\binits{P.~M.}}
(\byear{2000}).
\btitle{Whittle pseudo-maximum likelihood estimation for nonstationary time
series}.
\bjournal{J. Amer. Statist. Assoc.}
\bvolume{95}
\bpages{1229--1243}.
\bid{issn={0162-1459}, mr={1804246}}
\bptok{imsref}%
\end{barticle}
%
\endbibitem

%b26 #&#
\bibitem{walker}
%
\begin{barticle}[mr]
\bauthor{\bsnm{Walker},~\bfnm{A.~M.}\binits{A.~M.}}
(\byear{1964}).
\btitle{Asymptotic properties of least-squares estimates of parameters
of the
spectrum of a stationary non-deterministic time-series}.
\bjournal{J. Austral. Math. Soc.}
\bvolume{4}
\bpages{363--384}.
\bid{issn={0263-6115}, mr={0171345}}
\bptok{imsref}%
\end{barticle}
%
\endbibitem

%b27 #&#
\bibitem{zygmund}
%
\begin{bbook}[mr]
\bauthor{\bsnm{Zygmund},~\bfnm{A.}\binits{A.}}
(\byear{1977}).
\btitle{Trigonometric Series}. %{V}ol. {I}, {II}}.
\bpublisher{Cambridge Univ. Press}, \baddress{Cambridge}.
%Volumes I and
%II bound together}.
\bptok{imsref}%
\end{bbook}
%
\endbibitem

\end{thebibliography}
\end{document}